\documentclass[12pt]{amsart}

\usepackage{breakurl}
\usepackage{graphicx}
\usepackage{amsmath}
\usepackage{amscd}
\usepackage{amsfonts}
\usepackage{amssymb}
\usepackage{color}
\usepackage{fullpage}
\usepackage{mathtools}
\usepackage[pagebackref,hypertexnames=false, colorlinks, citecolor=red,linkcolor=blue, urlcolor=red]{hyperref}
\usepackage{todonotes}
\usepackage{enumitem}

\numberwithin{equation}{section}

\newtheorem{theorem}{Theorem}[section]
\newtheorem{lemma}[theorem]{Lemma}
\newtheorem{proposition}[theorem]{Proposition}

\newtheorem{corollary}[theorem]{Corollary}

\theoremstyle{definition}
\newtheorem{example}[theorem]{Example}
\newtheorem{remark}[theorem]{Remark}

\newtheorem{problem}[theorem]{Problem}

\newcommand{\be}{\begin{equation}}
\newcommand{\ee}{\end{equation}}
\newcommand{\bes}{\begin{equation*}}
\newcommand{\ees}{\end{equation*}}

\newcommand{\cA}{\mathcal{A}}
\newcommand{\cB}{\mathcal{B}}

\newcommand{\cD}{\mathcal{D}}

\newcommand{\cH}{\mathcal{H}}

\newcommand{\cK}{\mathcal{K}}
\newcommand{\cL}{\mathcal{L}}

\newcommand{\cU}{\mathcal{U}}

\newcommand{\cW}{\mathcal{W}}

\newcommand{\bC}{\mathbb{C}}
\newcommand{\bD}{\mathbb{D}}

\newcommand{\bF}{\mathbb{F}}
\newcommand{\bM}{\mathbb{M}}
\newcommand{\bN}{\mathbb{N}}

\newcommand{\bR}{\mathbb{R}}
\newcommand{\bT}{\mathbb{T}}
\newcommand{\bZ}{\mathbb{Z}}

\newcommand{\ol}{\overline}

\newcommand{\re}{\operatorname{Re}}
\newcommand{\im}{\operatorname{Im}}

\newcommand{\ru}{{\mathrm{u}}}

\newcommand{\Aut}{\operatorname{Aut}}

\newcommand{\conv}{\operatorname{conv}}

\newcommand{\id}{\operatorname{id}}

\newcommand{\UCP}{\operatorname{UCP}}

\newcommand{\distance}{\operatorname{d}}

\begin{document}

\title{Dilations of unitary tuples}

 \author{Malte Gerhold}
 \address{M.G., Institut f\"ur Mathematik und Informatik \\
 Universit\"at Greifswald\\Walther-Rathenau-Stra\ss{}e 47 \\
 17489 Greifswald \\ Germany\\}
 \email{mgerhold@uni-greifswald.de}
\urladdr{\href{https://math-inf.uni-greifswald.de/institut/ueber-uns/mitarbeitende/gerhold/}{\url{https://math-inf.uni-greifswald.de/institut/ueber-uns/mitarbeitende/gerhold/}}}

 \author{Satish K. Pandey}
 \address{S.P., Faculty of Mathematics\\
 Technion - Israel Institute of Technology\\
 Haifa\; 3200003\\
 Israel}
 \email{satishpandey@campus.technion.ac.il}
\urladdr{\href{http://noncommutative.space/}
{\url{http://noncommutative.space/}}}

 \author{Orr Moshe Shalit}
 \address{O.S., Faculty of Mathematics\\
 Technion - Israel Institute of Technology\\
 Haifa\; 3200003\\
 Israel}
 \email{oshalit@technion.ac.il}
 \urladdr{\href{https://oshalit.net.technion.ac.il/}{\url{https://oshalit.net.technion.ac.il/}}}

 \author{Baruch Solel}
 \address{B.S., Faculty of Mathematics\\
 Technion - Israel Institute of Technology\\
 Haifa\; 3200003\\
 Israel}
 \email{mabaruch@technion.ac.il}

 \thanks{The work of M. Gerhold is partially supported by the DFG, project no.\  397960675.}
 \thanks{The work of S. K. Pandey is supported in part at the Technion by a fellowship of the Israel Council for Higher Education.}
 \thanks{The work of O.M. Shalit is partially supported by ISF Grants no.\ 195/16 and 431/20.
 }
 \subjclass[2010]{47A13, 46L54}
 \keywords{Dilations, free unitaries, matrix range, noncommutative tori}

 \addcontentsline{toc}{section}{Abstract}

\begin{abstract}
We study the space of all $d$-tuples of unitaries $u = (u_1, \ldots, u_d)$ using dilation theory and matrix ranges. 
Given two such $d$-tuples $u$ and $v$ generating, respectively, C*-algebras $\cA$ and $\cB$, we seek the minimal dilation constant $c = c(u,v)$ such that $u \prec cv$, by which we mean that there exist faithful $*$-representations $\pi \colon \cA \to B(\cH)$ and $\rho\colon \cB \to B(\cK)$, with $\cH\subseteq \cK$, such that for all $i$, $\pi(u_i)$ is equal to the compression $P_\cH \rho(cv_i)\big|_\cH$ of $\rho(cv_i)$ to $\cH$. 
This gives rise to a metric 
\[
\distance_{\mathrm D}(u,v) = \log \max\{c(u,v), c(v,u)\}
\] 
on the set of equivalence classes of $*$-isomorphic tuples of unitaries. 
We compare this metric to the metric $\distance_{\mathrm{HR}}$ determined by 
\[
\distance_{\mathrm{HR}}(u,v) = \inf\left\{\|u' - v'\| : u', v' \in B(\cH)^d, u' \sim u \textrm{ and } v' \sim v\right\}, 
\]
and we show the inequality
\[
\distance_{\mathrm{HR}}(u,v) \leq K \distance_{\mathrm{D}}(u,v)^{1/2} 
\] where $1/2$ is optimal.
When restricting attention to unitary tuples whose matrix range contains a $\delta$-neighborhood of the origin, then $\distance_{\mathrm{D}}(u,v) \leq d \delta^{-1} \distance_{\mathrm{HR}}(u,v)$, so these metrics are equivalent on the set of tuples whose matrix range contains some neighborhood of the origin. 
Moreover, these two metrics are equivalent to the Hausdorff distance between the matrix ranges of the tuples. 

For particular classes of unitary tuples we find explicit bounds for the dilation constant. 
For example, if for a real antisymmetric $d \times d$ matrix $\Theta = (\theta_{k,\ell})$ we let $u_\Theta$ be the universal unitary tuple $(u_1, \ldots, u_d)$ satisfying $u_\ell u_k = e^{i\theta_{k,\ell}} u_k u_\ell$, then we find that $c(u_\Theta, u_{\Theta'}) \leq e^{\frac{1}{4}\|\Theta - \Theta'\|}$. 
Combined with the above equivalence of metrics, this allows to recover
the result of Haagerup-R{\o}rdam (in the $d=2$ case) and Gao (in the $d\geq 2$ case) that there exists a map $\Theta \mapsto U(\Theta) \in B(\cH)^d$ such that $U(\Theta) \sim u_\Theta$ and 
\[
\|U(\Theta) - U({\Theta'})\| \leq K \|\Theta - \Theta'\|^{1/2} .
\]

Of special interest are: the universal $d$-tuple of noncommuting unitaries $\ru$, the $d$-tuple of free Haar unitaries $u_f$, and the universal $d$-tuple of commuting unitaries $u_0$. 
We find upper and lower bounds on the dilation constants among these three tuples, 
and in particular we obtain rather tight (and surprising) bounds 
\[
2\sqrt{1 - \frac{1}{d}} \leq c(u_f,u_0) \leq 2 \sqrt{1-\frac{1}{2d}}.
\]
From this, we recover Passer's upper bound for the universal unitaries $c(\ru,u_0) \leq \sqrt{2d}$. 
In the case $d = 3$ we obtain the new lower bound $c(\ru,u_0) \geq 1.858$, which improves on the previously known lower bound $c(\ru,u_0) \geq \sqrt{3}$. 
\end{abstract}

\maketitle

\section{Introduction}\label{sec:introduction}

Let $A = (A_1, \ldots, A_d)$ be a $d$-tuple of operators on a Hilbert space $\cH$ and let $B = (B_1, \ldots, B_d)$ be a $d$-tuple of operators on a Hilbert space $\cK \supseteq \cH$. 
We say that $A$ is {\em the compression of $B$ to $\cH$} if 
\be\label{eq:dilationAB}
A = P_\cH B \big|_\cH. 
\ee
By \eqref{eq:dilationAB} we mean that $A_i =  P_\cH B_i \big|_\cH$ for all $i=1, \ldots, d$, where $P_\cH$ denotes the orthogonal projection $P_\cH \colon \cK \to \cH$.
In this case we say that $B$ is a {\em dilation} of $A$, and we write $A \prec B$.
Consider the following problem. 

\begin{problem}\label{prob:dilconst}
Fix $d \in \bN$.
What is the smallest constant $c$ such that for every $d$-tuple of contractions $A$, there exists a $d$-tuple of commuting normal contractions $B$ such that $A\prec cB$?
\end{problem}
Let us write $C_d$ for the smallest constant $c$ which is the solution to the above problem. 
Problem \ref{prob:dilconst} and similar problems have come up in the setting of relaxation of spectrahedral inclusion problems \cite{HKMS19}, in interpolation problems for completely positive maps and the study of the structure of operator systems \cite{DDSS17,FNT17}, and fit in the general paradigm of studying operator theory through dilations \cite{ShDilationSurvey}. 
More recently, such problems have turned out to be connected to quantum information theory \cite{BN18,BN20} as well as other aspects of mathematical physics and C*-algebras \cite{GS20}.

Passer, Shalit and Solel showed that if $A$ is a $d$-tuple of selfadjoint contractions, then there exists a $d$-tuple of commuting selfadjoint contractions $B$ such that $A \prec cB$ as in Problem \ref{prob:dilconst} with $c = \sqrt{d}$, and that this is the optimal constant for selfadjoint tuples \cite[Theorem 6.6]{PSS18}.
Moreover, it was shown by Passer in \cite[Theorem 4.4]{Passer} that if $A$ is not assumed selfadjoint, one can do with $c = \sqrt{2d}$.
Thus, we have the bounds
\be\label{eq:Cd}
\sqrt{d} \leq C_d \leq \sqrt{2d}.
\ee

\begin{remark}\label{rem:contraction_unitaries}
It is convenient to note that Problem \ref{prob:dilconst} can be reformulated in terms of unitaries as follows: {\em What is the smallest constant $c = C_d$ such that for every $d$-tuple of {\bf unitaries} $V$, there exists a $d$-tuple of commuting {\bf unitaries} $U$, such that $V\prec cU$?}

To see why, when replacing the tuple $A$ of contractions from Problem \ref{prob:dilconst} with a tuple of unitaries $V$, we end up with an equivalent problem, we note that the simple construction 
\[
V_i = 
\begin{pmatrix}
      A_i & (1 - A_i A_i^*)^{1/2}\\
      (1 - A_i^* A_i)^{1/2} & -A_i^* 
\end{pmatrix}
\]
provides a unitary dilation $A \prec V$ for any tuple of contractions $A$. 
Moreover, by a minor modification of the proof of \cite[Proposition 2.3]{PSS18}, if we can dilate to a tuple of commuting normals $B$ with joint spectrum contained in a compact convex set $K$, then we can dilate to a commuting normals with joint spectrum contained in $\overline{\operatorname{ext}(K)}$. 
Thus, if there a is dilation $cB$ where $B$ is a tuple of commuting normal contractions (so that $\sigma(B) \subseteq \ol{\bD}^d$), then there is also a dilation of the form $cU$, where $U$ is a tuple of commuting unitaries. 
\end{remark}

The goal of this paper is to study dilation theory in the context of unitary tuples with Problem \ref{prob:dilconst} in sight. 
The next subsection is devoted to setting the notation for the rest of the paper, and in the following one we will give a summary of our main results.

\subsection*{Some definitions and notation}

In this paper, $d$ will always be some positive integer that may be considered as fixed throughout. If not indicated otherwise, sums will be assumed to run from $1$ to $d$ over all appearing indices. 
Our main concern will be $d$-tuples of unitaries. 
For a Hilbert space $\cH$, we let $B(\cH)$ denote the algebra of bounded operators on $\cH$, and $\cU(\cH)$ the group of unitary operators on $\cH$. 
Likewise $B(\cH)^d$ and $\cU(\cH)^d$ will denote $d$-tuples of operators or unitaries, respectively, on $\cH$.
We let $M_n = M_n(\bC)$ denote the set of all $n \times n$ matrices over $\bC$, and $M_n^d$ the set of all $d$-tuples of such matrices.
The ``noncommutative universe'' (in $d$ variables) is the disjoint union $\bM^d = \bigcup_{n=1}^\infty M_n^d$. 
The {\em matrix range} \cite{Arv72} of a tuple $A = (A_1, \ldots, A_d)$ in $B(\cH)^d$ is the disjoint union $\cW(A) = \bigcup_{n} \cW_n(A)$, where for all $n\in\mathbb N$, the set $\cW_n(A) \subseteq M_n^d$ is defined by
\[
\cW_n(A) = \bigl\{\bigl(\phi(A_1), \ldots, \phi(A_d)\bigr) : \phi \in \UCP\bigl(B(\cH), M_n\bigr)\bigr\};
\]
here and below, UCP stands for {\em unital completely positive}, and $\UCP(B(\cH), M_n)$ is the set of all UCP maps from $B(\cH)$ to $M_n$.

For every $d$-tuple of operators $X = (X_1, \ldots, X_d)$ we write $\|X\| := \max_i \left\| X_i \right\|$.
The space $\cB(\cH)^d$ then becomes a metric space with $\distance(X,Y) = \left\|X-Y\right\|$. 
The norm on $M_n^d$ also induces a distance function on the subsets of $M_n^d$, the {\em Hausdorff distance}, given by
\[
\distance_{\mathrm H}(E,F) = \max\left\{\sup_{x \in E} \distance(x,F),  \sup_{y \in F}  \distance(y,E)\right\},
\]
where $\distance(x,F) = \inf_{y \in F}\distance(x,y)$. 
The {\em Hausdorff metric}, also denoted $\distance_{\mathrm H}$, is the restriction of the Hausdorff distance to the compact subsets of $M_n^d$.  

We denote by $\mathbb{F}_d$ the free group on $d$ generators and by $C_r^*(\bF_d)$ its reduced C*-algebra. 
We let $\ru = (\ru_1, \ldots, \ru_d)$ denote the universal $d$-tuple of noncommuting unitaries generating the full C*-algebra $C^*(\bF_d)$ of the free group. 
We let $u_f = (u_{f,1}, \ldots, u_{f,d})$ denote the canonical $d$-tuple of unitaries generating $C_r^*(\bF_d)$, and let $\tau_f$ be the canonical tracial state on $C_r^*(\bF_d)$. 
We shall refer to any $d$-tuple $u$ of unitaries generating a tracial C*-algebra $(\cA, \tau)$ such that $u_f$ and $u$ have the same $*$-distributions with respect to $\tau_f$ and $\tau$ as \emph{free Haar unitaries}.
Finally, we write $u_0 = (u_{0,1}, \ldots, u_{0,d})$ for the universal commuting unitary $d$-tuple (i.e., the canonical generators of the commutative C*-algebra $C^*(\bZ^d) \cong C(\bT^d)$). 

Whenever we shall want to explicitly refer to the ``number of variables'' $d$ being considered, we shall write $\ru^{(d)}$, $u_f^{(d)}$, and $u_0^{(d)}$, etc., to emphasize this.

\subsection*{Summary of the main results} 

In the next section we consider three distance functions on the set of (equivalence classes of) $d$-tuples of unitaries: 
\[
\distance_{\mathrm{HR}}(u,v):=\inf\left\{\left\|u'-v'\right\|\colon u', v' \in B(\cH)^d,  u \sim u' \textrm{ and } v \sim v'\right\}, 
\]
\[
\distance_{\mathrm{mr}}(u,v):=\distance_{\mathrm H}(\cW(u),\cW(v)) , 
\]
and 
\[
\distance_{\mathrm{D}}(u,v):= \log \max\{c(u,v), c(v,u)\} , 
\]
where $c(u,v) = \inf\{c : u \prec cv\}$. 
As we explained in Remark \ref{rem:contraction_unitaries}, $c(\ru,u_0)$ is equal to $C_d$. 
Although we have not succeeded in determining $C_d$, we argue that understanding $c(u,v)$ for various $d$-tuples of unitaries is beneficial for finding $C_d$. 
Hopefully, the readers will soon be convinced that studying the constants $c(u,v)$ has interesting consequences which a priori seem to have nothing to do with dilation theory. 

We prove that the above distance functions are metrics, and that they are equivalent on tuples whose matrix range contains a neighborhood of the origin. 
In Theorem \ref{thm:dilandHR} we show that there exists a constant $K$ such that
\[
\distance_{\mathrm{HR}}(u,v) \leq K \distance_{\mathrm{D}}(u,v)^{1/2}. 
\]
In Corollary \ref{cor:metricsU_delta} we show that $\distance_{\mathrm{D}}(u,v) \leq d \delta^{-1} \distance_{\mathrm{HR}}(u,v)$ holds, when restricting attention to unitary tuples whose matrix range contains a $\delta$-neighborhood of the origin. 

In Section \ref{sec:free} we focus on dilation constants involving the free Haar unitaries. 
Thanks to the highly developed machinery of free probability, we are able to obtain some exact values and close estimates. 
In the direction of dilating to free unitaries, we obtain the exact values 
\[
c(\mathrm{u},u_f) = c(u_0,u_f) = \frac{d}{\sqrt{2d-1}} \,\,,
\]
see Corollary \ref{cor:c_commUni_free}. 
In the other direction, i.e., dilating free unitaries to commuting unitaries, we obtain the bounds 
\[
2\sqrt{1-\frac{1}{d}} \leq c(u_f, u_0) \leq 2\sqrt{1-\frac{1}{2d}} \,\,,
\] 
see Theorems \ref{thm:cf0-geq-sqrt2} and \ref{thm:low_freecomm}. 
The combination of Corollary \ref{cor:c_commUni_free} and Theorem \ref{thm:low_freecomm} allows us, in Corollary \ref{cor:Passer}, to recover Passer's bound $C_d = c(\ru^{(d)},u^{(d)}_0) \leq \sqrt{2d}$ (see \cite[Theorem 4.4]{Passer}). 

Another interesting class of unitary tuples are the noncommutative tori, to which we turn in Sections \ref{sec:NCT_abs} and \ref{sec:NCT_cont}. 
Given an antisymmetric real matrix $\Theta=(\theta_{k,\ell})$, let $u_\Theta = (u_{\Theta,1}, \ldots, u_{\Theta,d})$ denote the universal unitary $d$-tuple that satisfies the commutation relations
\[
u_{\Theta,\ell} u_{\Theta,k}=e^{i\theta_{k,\ell}}u_{\Theta,k} u_{\Theta,\ell} \,\, , \,\, k,\ell = 1, \ldots, d .
\]
Note that $u_0$ --- the universal $d$-tuple of commuting unitaries --- corresponds to the noncommutative torus with $\Theta = 0_{d\times d}$. 
In Section \ref{sec:NCT_abs} we collect some results about the constant $c(u_\Theta,u_{\Theta'})$. 
For example, in Proposition \ref{prop:sametheta} we find that when $\theta_{k,\ell} = \theta$ for all $k<\ell$, we have the value
\[
c(u_\Theta,u_0) = \frac{2d}{\|u_{\Theta,1}+u_{\Theta,1}^* + \ldots + u_{\Theta,d} + u_{\Theta,d}^*\|}.
\]
This result generalizes \cite[Theorem 6.3]{GS20}, where the same formula was obtained for $d = 2$. 
In \cite{GS20} this formula was used with explicit computations to show that 
\[
C_2 \geq \sup_\theta c((u_{\theta,1},u_{\theta,2}), u_0) \geq 1.543 > \sqrt{2} \,\,, 
\]
thereby showing that the lower bound in \eqref{eq:Cd} is not sharp, at least for $d = 2$. 
Similarly, in Corollary \ref{cor:lower_bound} we use the above formula as the theoretical backbone of some numerical computations which show that
\[
C_3 \geq \sup_\Theta c(u_\Theta, u_0) \geq 1.858 > \sqrt{3} \,\,,
\]
which shows the same for $d=3$. 

In Section \ref{sec:NCT_cont} we follow \cite[Section 3]{GS20} and find representations of the noncommutative tori using the Weyl unitaries. 
These representations are used to show that $c(u_{\Theta},u_{\Theta'}) \leq e^{\frac{1}{4}\|\Theta-\Theta'\|}$ (Theorem \ref{thm:Weyl-rep}) and, consequently, $\distance_{\mathrm{D}}(u_{\Theta},u_{\Theta'}) \leq \frac{1}{4}\|\Theta-\Theta'\|$. 
In combination with Theorem \ref{thm:dilandHR}, this yields (Corollary \ref{cor:nctori_cont}) the existence of a constant $K$ such that 
\[
\distance_{\mathrm{HR}}(u_\Theta,u_{\Theta'}) \leq K \|\Theta-\Theta'\|^{1/2} .
\] 
Thus, there is a continuous path $\Theta \mapsto u_\Theta$ from the space of real antisymmetric matrices into the metric space determined by $\distance_{\mathrm{HR}}$ which is H{\"o}lder continuous with exponent $1/2$. 
In the appendix to this paper, we show that the techniques of Haagerup and R{\o}rdam from \cite{HR95} can be used to prove that every H{\"o}lder continuous path $t \mapsto u_t$ from the interval $[0,1]$ into the space of equivalence classes of unitaries, endowed with the metric $\distance_{\mathrm{HR}}$, can be lifted to a path of representations $t \mapsto U(t) \in \cU(\cH)^d$ on the same Hilbert space $\cH$ which is H{\"o}lder continuous (with the same exponent) with respect to the operator norm. 
Together with Corollary \ref{cor:nctori_cont}, this allows us to recover the result from \cite{Gao18,HR95} that there exists a norm continuous map $\Theta \mapsto U(\Theta) \in B(\cH)^d$, such that $U(\Theta) \sim u_\Theta$ and
\[
\|U(\Theta) - U({\Theta'})\| \leq K \|\Theta - \Theta'\|^{1/2} .
\]

\section{Metric structure on the space of unitary tuples}\label{sec:metric} 

For a pair of unitary tuples $u = (u_1, \ldots, u_d)$ and $v =(v_1, \ldots, v_d)$, we say that $u$ is {\em equivalent} to $v$, and we write $u \sim v$, if there exist a $*$-isomorphism $\pi \colon C^*(u) \to C^*(v)$ such that $\pi(u_i) = v_i$ for all $i=1,\ldots, d$.
By Voiculescu's theorem (see, e.g., \cite[Corollary II.5.6]{DavBook}), we have that $u\sim v$ if and only if the infinite ampliations of $u$ and $v$ are approximately unitarily equivalent, a fact that will be important in the appendix of this paper. 
For the main part of the paper, we will make extensive use of a characterization of equivalence in terms of matrix ranges.  
Note that equivalent unitary tuples have the same matrix ranges. In fact, we will show soon that the converse also holds, so that unitary tuples are equivalent if and only if their matrix ranges agree (see Proposition \ref{prop:metrics}). 

Let $\cU(d)$ be the set of all equivalence classes of $d$-tuples of unitaries. 
It suffices to consider only separably acting unitaries, because every separable C*-algebra can be represented faithfully on a separable Hilbert space. 
Let $\cU_\delta(d)$ be the subset of $\cU(d)$ of all unitary $d$-tuples $u$ such that $\cW_1(u)$ contains the neighborhood $B_\delta(0) = \{z \in \bC^d : \|z\|<\delta\}$ of the origin, and define $\cU_0(d)$ as the subset of $\cU(d)$ that consists of all tuples $u$ such that $\cW_1(u)$ contains some neighborhood of the origin. 
We shall usually identify tuples with the equivalence classes that they belong to, unless there is a special reason to be careful. 

\begin{example}\label{expl:gauge}
We wish to show that $\cU_0(d)$ contains a rich class of interesting examples. 
Suppose that $u \in \cU(d)$ is such that $C^*(u)$ carries a natural \emph{gauge action}, that is, there is homomorphism $\gamma \colon \bT^d \to \Aut (C^*(u))$ such that $\gamma_\lambda(u_i) = \lambda_i u_i$ for every $\lambda = (\lambda_1, \ldots, \lambda_d) \in \bT^d$ and all $i=1, \ldots, d$. 
Many unitary tuples of interest carry such an action, for example the canonical generators of the full and reduced free group C*-algebras, the generators of $C(\bT^d)$, and the canonical generators of the noncommutative tori. 
Since for every $i$, there is a state $\phi$ of $C^*(u)$ such that $|\phi(u_i)| = 1$, we can integrate away the other coordinates to find a state $\psi$ such that $|\psi(u_j)| = \delta_{ij}$. 
By gauge invariance and convexity of $\cW(u)$, we see that $\cW_1(u)$ contains the unit ball of $\bC^d$ with respect to the $\ell^1$ norm. 
In conclusion: whenever $u$ has a gauge action as above, $\cW_1(u)$ contains $B_{1/\sqrt{d}}(0)$, and thus $u \in \cU_{1/\sqrt{d}}(0)$. 
\end{example}

Given two unitary tuples $u,v$ and a positive real number $c$, we write $u \prec cv$ if there exist two Hilbert spaces $\cH \subseteq \cK$ and two operator tuples $U \in B(\cH)^d$ and $V \in B(\cK)^d$, such that $u \sim U$, $v \sim V$ and 
\[
U = P_\cH cV \big|_\cH. 
\]
By Stinespring's theorem, $u \prec cv$ if and only if there is a UCP map from the operator system generated by $v$ to the operator system generated by $u$, that maps $cv$ to $u$. 

For $u,v \in \cU(d)$, we put 
\[
c(u,v) = \inf\{c : u \prec cv\}. 
\]
This infimum is actually attained, by compactness of UCP maps in the BW topology (see \cite[Theorem 7.4]{PauBook}). 
By Remark \ref{rem:contraction_unitaries}, if $\ru$ denotes the universal noncommuting unitary $d$-tuple (i.e., the canonical generators of the full C*-algebra of the free group $C^*(\bF_d)$), and if $u_0$ denotes the universal commuting unitary $d$-tuple (i.e., the canonical generators of the commutative C*-algebra $C^*(\bZ^d) \cong C(\bT^d)$), then $c(\ru,u_0) = C_d$. 

We shall consider the following distance functions on $\cU(d)$. 
\begin{description}
\item[The Haagerup-R{\o}rdam-distance]
\[
\distance_{\mathrm{HR}}(u,v):=\inf\left\{\left\|u'-v'\right\|\colon u', v' \in B(\cH)^d,  u \sim u' \textrm{ and } v \sim v'\right\}. 
\]
\item[The dilation distance]
\[
\distance_{\mathrm{D}}(u,v):= \log \max\bigl\{c(u,v), c(v,u)\bigr\} .
\]
\item[The Matrix Range distance]
\[
\distance_{\mathrm{mr}}(u,v):=\distance_{\mathrm H}\bigl(\cW(u),\cW(v)\bigr) .
\]
\end{description}

Our goal in this section is to understand these distances and the relationships between them. 
We will show that they are all metrics. 
The metric $\distance_{\mathrm{HR}}$ is inspired by the metric introduced by Haagerup and R{\o}rdam in \cite[Definition 4.1]{HR95}, and hence the terminology. Note that, trivially, $\distance_{\mathrm{HR}}(u,v)\leq 2$ for all $u,v\in \cU(d)$.

\begin{lemma}\label{lem:metrics}
$\distance_{\mathrm{mr}} \leq \distance_{\mathrm{HR}}$ and $\distance_{\mathrm{mr}} 
\leq 2\distance_{\mathrm D}$. 
\end{lemma}

\begin{proof}
The first inequality follows readily from the facts that the matrix ranges of equivalent unitary tuples agree and that UCP maps are contractive. 
The second one follows from the definitions by a short calculation: if $u \prec cv$, and $X \in \cW(u) \subseteq c\cW(v)$, then $X = cY$ for $Y \in \cW(v)$, and $\|X - Y \| = (c-1)\|Y\| \leq c-1$. 
This argument shows that $\distance_{\mathrm H}(\cW(u),\cW(v)) = \distance_{\mathrm{mr}}(u,v) \leq \max\{c(u,v),c(v,u)\} - 1 = e^{\distance_{\mathrm{D}}(u,v)}-1$ (see also \cite[Proposition 2.4]{PP+}). Finally, note that $e^x-1\leq 2x$ for all $x\in[0,1]$, so the last claimed inequality holds for all $u,v$ with $\distance_{\mathrm D}(u,v)\leq 1$. On the other hand, if $d_D(u,v)>1$, then
$\distance_{\mathrm{mr}}(u,v)\leq \distance_{\mathrm{HR}}(u,v)\leq 2 < 2\distance_{\mathrm D}(u,v)$ holds trivially.
\end{proof}

\begin{proposition}\label{prop:metrics}
The distance functions $\distance_{\mathrm{HR}}$, $\distance_{\mathrm{D}}$, and $\distance_{\mathrm{mr}}$ are all metrics on $\cU(d)$. 
\end{proposition} 

\begin{proof}
For each of the three distance functions, the only nontrivial part is to show that $d(u,v) = 0 \implies u \sim v$. 
As $\distance_{\mathrm{mr}} \leq \distance_{\mathrm{HR}}$ and $\distance_{\mathrm{mr}} \leq 2\distance_{\mathrm{D}}$, it suffices to consider $\distance_{\mathrm{mr}}$. 

Thus, suppose that $\distance_{\mathrm{mr}}(u,v) = 0$. 
By \cite[Theorem 5.1]{DDSS17}, there exists a unital and completely isometric map $\phi\colon S_u \to S_v$ from the operator system $S_u$ generated by $u$ to the operator system $S_v$ generated by $v$ such that $\phi(u) = v$. 
By \cite[Corollary 2.2.8]{Arv69}, the Shilov ideal of an operator system generated by unitaries, relative to the C*-algebra that the unitaries generate, is trivial. 
Thus, by \cite[Theorem 2.2.5]{Arv69}, $\phi$ is equal to the restriction of a C*-isomorphism $\pi \colon C^*(u) \to C^*(v)$. 
We conclude that $v = \pi(u)$, that is, $v \sim u$. 
\end{proof}

By Lemma \ref{lem:metrics}, we see that if a sequence $\{u^{(n)}\}$ in $\cU(d)$ converges to $u$ with respect to $\distance_{\mathrm{HR}}$ or $\distance_{\mathrm{D}}$, then it also converges with respect to $\distance_{\mathrm{mr}}$. 
We will now see that the converse holds in $\cU_0(d)$. 
In fact, we will see that for every $\delta > 0$ the three metrics are strongly equivalent in $\cU_\delta(d)$. 

\begin{proposition}\label{prop:metricsU_delta}
Let $\delta > 0$. 
Then for all $u,v \in \cU_\delta(d)$, 
\[
\distance_{\mathrm{D}}(u,v) \leq d \delta^{-1} \distance_{\mathrm{mr}}(u,v). 
\]
\end{proposition}

\begin{proof}
Let us assume that $C^*(u)$ is faithfully represented on the Hilbert space $\cH$. 
By a simple geometric argument (see \cite[Proposition 2.4]{PP+}), one shows that 
\[
\cW(u) \subseteq c \cW(v)
\]
for 
\[
c = 1+ d \delta^{-1} \distance_{\mathrm{mr}}(u,v). 
\] 
Thus, there is a UCP map $\phi \colon C^*(v) \to B(\cH)$ such that $\phi(cv) = u$. 
By Stinespring's theorem, there is a $*$-representation $\pi \colon C^*(v) \to B(\cK)$, where $\cK \supseteq \cH$, such that $P_\cH c\pi(v)\big|_\cH = u$. 
We may assume (by adding to $\pi$ a direct summand which is faithful) that $\pi$ is faithful, and thus $u \prec cv$. 
The argument works with the roles of $u$ and $v$ reversed, thus
\[
\distance_{\mathrm D}(u,v)=\log\max\bigl\{c(u,v), c(v,u)\bigr\} \leq \log\bigl(1+ d \delta^{-1} \distance_{\mathrm{mr}}(u,v)\bigr)\leq d \delta^{-1} \distance_{\mathrm{mr}}(u,v)\bigr).\qedhere
\] 
\end{proof}

Recall that two metrics on the same space are called \emph{equivalent} if they generate the same topology and \emph{strongly equivalent} if they dominate each other up to constants. By Lemma \ref{lem:metrics} and Proposition \ref{prop:metricsU_delta} we get the following immediate corollary. 

\begin{corollary}\label{cor:metricsU_delta}
For every $\delta > 0$, the restrictions of $\distance_{\mathrm{D}}$ and $\distance_{\mathrm{mr}}$ to $\cU_\delta(d)$ are strongly equivalent metrics. Consequently, their restrictions to $\cU_0(d)$ are equivalent metrics.
\end{corollary} 

Now we wish to show that the dilation distance dominates the Haagerup-R{\o}rdam distance. 

\begin{theorem}\label{thm:dilandHR}
  If $u\prec cv$ and $v\prec cu$ then $\distance_{\mathrm{HR}}(u,v)\leq 28\sqrt{1-c^{-2}}$
  and, consequently,
  \[\distance_{\mathrm{HR}}(u,v)\leq 56 \distance_{\mathrm{D}}(u,v)^{1/2}.\]
\end{theorem} 

\begin{proof}
We assume
\be\label{eq:equiv}
cu \sim
\begin{pmatrix}
v & x \\ y & z 
\end{pmatrix} \quad \textrm{and} \quad
cv \sim 
\begin{pmatrix}
u & r \\ s & t 
\end{pmatrix}.
\ee
As $u$ and $v$ are unitary, we find that $\|x\| = \|y\| = \|s\| = \|t\| = \sqrt{c^2 - 1}$. 
We shall keep this in mind, as in the rest of the proof we will be able to bound off-diagonal block by some constant times $\sqrt{c^2 - 1}$

We will need to be careful and track the identifications made. 
But first, we replace every operator $a$ appearing above with the infinite ampliation $a \oplus a \oplus \cdots$, so we may assume that $u,v,x,y$, etc., are all given as concrete operators acting on an infinite dimensional Hilbert space $\cH$. 
So the equivalences in \eqref{eq:equiv} are due to $*$-isomorphisms $\pi \colon C^*(u) \subset B(\cH) \to B(\cH \oplus \cH)$ and $\rho \colon C^*(v) \subset B(\cH) \to B(\cH \oplus \cH)$ such that 
\[
c\pi(u) =
\begin{pmatrix}
v & x \\ y & z 
\end{pmatrix} \quad \textrm{and} \quad
c\rho(v) = 
\begin{pmatrix}
u & r \\ s & t 
\end{pmatrix}.
\]
In fact, by applying the standard combination of Arveson's extension theorem and Stinespring's dilation theorem, we may assume that $\pi$ and $\rho$ are $*$-homomorphisms defined on all of $B(\cH)$. 
Letting $\pi^{(k)}$ and $\rho^{(k)}$ denote the ampliations of the representations, we obtain
\begin{align*}
c^2 \pi^{(2)}\rho(v) 
&= c\begin{pmatrix} \pi(u) & \pi(r) \\ \pi(s) & \pi(t)  \end{pmatrix} \\
&= \begin{pmatrix} \begin{bmatrix}v & x \\ y & z \end{bmatrix} & c\pi(r) \\ c\pi(s) & c\pi(t) \end{pmatrix} \in B(\cH^4) . 
\end{align*}
On the other hand, we find that 
\begin{align*}
c^3 \pi^{(4)} \rho^{(2)} \pi(u)
&= c^2 \pi^{(4)} \rho^{(2)} \begin{pmatrix} v & x \\ y & z \end{pmatrix} \\ 
&= c \pi^{(4)} \begin{pmatrix} \begin{bmatrix}u & r \\ s & t \end{bmatrix} & c\rho(x) \\ c\rho(y) & c \rho(z) \end{pmatrix} \\
&= \begin{pmatrix} \begin{bmatrix} \begin{bmatrix} v & x \\ y & z \end{bmatrix} & c \pi(r) \\ c \pi(s) & c \pi(t) \end{bmatrix} & c^2 \pi^{(2)} \rho(x) \\ c^2 \pi^{(2)} \rho(y) & c^2 \pi^{(2)} \rho(z) \end{pmatrix} \in B(\cH^8) . 
\end{align*}
To summarize a little more concisely what we found: 
\[
v \sim  \begin{pmatrix} \frac{1}{c^2}v & * & * \\ * & \frac{1}{c^2}z & * \\ * & * & \frac{1}{c} \pi(t) \end{pmatrix} \in B(\cH \oplus \cH \oplus \cH^2), 
\]
and
\[
u \sim  \begin{pmatrix} \frac{1}{c^3}v & * & * & *  \\ * & \frac{1}{c^3}z & * & * \\ * & * & \frac{1}{c^2} \pi(t) & * \\ * & * & * & \frac{1}{c}\pi^{(2)} \rho(z) \end{pmatrix} \in B(\cH \oplus \cH \oplus \cH^2 \oplus \cH^4) , 
\]
where all off diagonal blocks denoted by $*$ are of norm less than $\frac{1}{c} \sqrt{c^2-1}$. 
Applying a permutation on the direct summands, it is convenient to rewrite this as follows: 
\[
v \sim V := \begin{pmatrix} \frac{1}{c^{2}}v & * & * \\ * & \frac{1}{c}\pi(t) & * \\ * & * & \frac{1}{c^{2}}z \end{pmatrix} \in B(\cH  \oplus \cH^2 \oplus \cH ), 
\]
and
\[
u \sim U: =  \begin{pmatrix} \frac{1}{c^3}v & * & * & *  \\ * & \frac{1}{c^2} \pi(t) & * & * \\ * & * & \frac{1}{c^3}z & * \\ * & * & * & \frac{1}{c}\pi^{(2)} \rho(z) \end{pmatrix} \in B(\cH \oplus \cH^2  \oplus \cH \oplus \cH^4) .
\]
If we restrict $\pi^{(2)} \rho$ to the C*-algebra generated by $z$ (which we have assumed to have infinite multiplicity), then we are in the situation of Voiculescu's theorem, which tells us that the representations $\id$ and ${\id} \oplus \pi^{(2)} \rho$ are approximately unitarily equivalent, that is, there is a sequence of unitaries $w_n \colon \cH \to \cH \oplus \cH^4$ such that $\lim_{n\to \infty} \|w_n z w_n^* - z \oplus \pi^{(2)} \rho(z) \| = 0$ (see \cite[Corollary II.5.5]{DavBook}). 
Then letting $W_n = I_{\cH} \oplus I_{\cH^2} \oplus w_n$, we have 
\[
\limsup_{n\to \infty} \left\| W_n V W_n^*- U \right\| \leq  \frac{24}{c}\sqrt{c^2-1} + \left(\frac{1}{c^2} - \frac{1}{c^3} \right)(2+c) + \left(\frac{1}{c} - \frac{1}{c^2} \right), 
\]
where the term $\frac{24}{c}\sqrt{c^2-1}$ accounts for the off diagonal blocks, the term with $2+c$ accounts for the first three diagonal blocks, and the term $\left(\frac{1}{c} - \frac{1}{c^2} \right)$ corresponds to the difference in the lower right block. 
Using $c \geq 1$ in order to simplify, we have $c-1\leq \sqrt{c^2-1}$ and $\frac{2+c}{c}\leq 3$, so doing the math with generous bounds, we find that 
\[
\limsup_{n\to \infty} \| W_n V W_n^*- U \|\leq \frac{28}{c}\sqrt{c^2-1}=28 \sqrt{1-c^{-2}} .
\]
Since $U \sim u$ and $W_n V W_n^* \sim v$ for all $n$, we conclude that 
\[
\distance_{\mathrm{HR}}(u,v) \leq 28\sqrt{1-c^{-2}}, 
\]
as required.

Now, letting $\delta = \distance_{\mathrm{D}}(u,v)$, we find that $\distance_{\mathrm{HR}}(u,v) \leq 28\sqrt{1-e^{-2\delta}}$.
Restricting attention to the bounded interval $\delta \in [0,1/2]$, we easily obtain 
\[
28\sqrt{1-e^{-2\delta}}\leq 28 \sqrt{e^{2\delta}-1} \leq 28\sqrt{4\delta} = 56 \distance_{\mathrm{D}}(u,v)^{1/2}.
\]
On the other hand, for $\delta>1/2$, the inequality $28\sqrt{1-e^{-2\delta}}\leq 28\leq 56 \distance_{\mathrm{D}}(u,v)^{1/2}$
is obvious. 
\end{proof}

\begin{remark}
Although we have not attempted to find the value of the optimal constant $K$ in the inequality $\distance_{\mathrm{HR}}(u,v)\leq K\distance_{\mathrm{D}}(u,v)^{1/2}$, it is interesting to note that the exponent $1/2$ is indeed optimal. 
To see this, we consider the universal pair of unitaries $u_\theta = (U_\theta,V_\theta)$ that satisfy $V_\theta U_\theta = e^{i\theta}U_\theta V_\theta$. 
By \cite[Theorem 3.2]{GS20}, we have $\distance_{\mathrm{D}}(u_\theta, u_{\theta'}) \leq 1/4|\theta - \theta'|$. 
On the other hand, by \cite[Proposition 4.6]{HR95}, $\distance_{\mathrm{HR}}(u_\theta, u_{\theta'}) \geq 1/2|\theta - \theta'|^{1/2}$; it follows that the inequality $\distance_{\mathrm{HR}}(u,v)\leq K\distance_{\mathrm{D}}(u,v)^{\alpha}$ cannot hold for any constant $K>0$ with $\alpha > 1/2$. 
\end{remark}

\begin{corollary}\label{cor:HRandmr}
The restrictions of $\distance_{\mathrm{HR}}$ and $\distance_{\mathrm{mr}}$ to $\cU_0(d)$ are equivalent. 
\end{corollary}

\begin{proof}
Lemma \ref{lem:metrics} recorded the trivial inequality $\distance_{\mathrm{mr}}\leq \distance_{\mathrm{HR}}$. 
On the other hand, if $u^{(n)}$ and $u$ are in $\cU_0(d)$ such that $\distance_{\mathrm{mr}}(u^{(n)},u) \rightarrow 0$, then $u \in \cU_\delta(d)$ for some $\delta > 0$, and therefore also $u^{(n)} \in \cU_\delta(d)$ for sufficiently large $n$. 
Combining Theorem \ref{thm:dilandHR} and Proposition \ref{prop:metricsU_delta}, we find that $\distance_{\mathrm{HR}}(u^{(n)},u) \rightarrow 0$, as required. 
\end{proof}

It is interesting to compare Corollary \ref{cor:HRandmr} to Theorems 2.3 and 2.4 in \cite{GS+}, which say that for a family of unitary tuples $u(t) = (u_1(t), \ldots, u_d(t))$, the {\em levelwise} convergence 
\[
\lim_{t\to t_0} \distance_{\mathrm H}(\cW_n(u(t)),\cW_n(u(t_0))) = 0 
\] 
for all $n$ is equivalent to $\lim_{t \to t_0} \|p(u(t))\| = \|p(u(t_0))\|$ for every $*$-polynomial $p$, i.e., the family $u(t)$ generates a continuous field of C*-algebras. 
On the other hand 
\[
\distance_{\mathrm{mr}}(u(t),u(t_0)) = \sup_{n \geq 1} \distance_{\mathrm H}(\cW_n(u(t)),\cW_n(u(t_0))), 
\]
and, by Corollary \ref{cor:HRandmr}, the uniform (over the levels) convergence $\lim_{t\to t_0} \distance_{\mathrm{mr}}(u(t),u(t_0)) = 0$ implies the stronger conclusion that appropriate $*$-isomorphic copies of $u(t)$ and $u(t_0)$ become as close in norm as we wish. 

\section{Dilations and the free Haar unitaries}\label{sec:free} 

\begin{lemma}\label{lem:moments}
Let $u = (u_1, \ldots, u_d)$ and $v = (v_1, \ldots, v_d)$ be two tuples of unitaries generating two C*-algebras $\cA$ and $\cB$, respectively. 
If $\tau$ and $\psi$ are faithful states on $\cA$ and $\cB$, respectively, and if $u$ and $v$ have the same $*$-distributions with respect to $\tau$ and $\psi$, then the map $u_i \mapsto v_i$ extends to a $*$-isomorphism between $\cA$ and $\cB$. 
\end{lemma} 

\begin{proof}
This is a well known fact in noncommutative probability. 
See, e.g., Theorem 2 on p.~163 in \cite{MSBook}. 
\end{proof}

As an immediate consequence, we obtain the following lemma. 
\begin{lemma}\label{lem:auto}
For every $\lambda = (\lambda_1, \ldots, \lambda_d) \in \bT^d$, and every $\pi$ in the symmetric group $S_d$, there exists an automorphism $\sigma \in \Aut (C_r^*(\bF_d))$ such that $\sigma(u_{f,i}) = \lambda_i u_{f,\pi(i)}$ for $i = 1, \ldots, d$. 
\end{lemma}
As in Example \ref{expl:gauge}, we call an automorphism $\sigma$ of the form $\sigma(u_{f,i}) = \lambda_i u_{f,i}$ ($i=1, \ldots, d$) a {\em gauge automorphism}, and we write $\sigma = \gamma_\lambda$. 
If $\sigma(u_{f,i}) =  u_{f,\pi(i)}$ for all $i$ (that is, if all $\lambda_i$s are equal to $1$), then we write $\sigma = \sigma_\pi$. 

A second useful consequence of Lemma \ref{lem:moments} is the following lemma. 
\begin{lemma}\label{lem:tensor_f}
Let $v = (v_1, \ldots, v_d)$ be a tuple of unitaries generating a C*-algebra $\cA$, $\varphi$ a faithful state on $\cA$, and put $u_f \otimes v = (u_{f,1} \otimes v_1, \ldots, u_{f,d} \otimes v_d)$. 
Then $\tau_f\otimes\varphi$ is a faithful trace on $C^*(u_f \otimes v)$, and there exists a trace-preserving $*$-isomorphism from $(C_r^*(\bF_d), \tau_f)$ onto $(C^*(u_f \otimes v), \tau_f \otimes \varphi)$ sending $u_{f,i}$ to $u_{f,i} \otimes v_i$. 
\end{lemma}

\begin{proof} 
The tensor product $\tau_f\otimes\varphi$ is a faithful trace on $C^*(u_f \otimes v)$ as the restriction of the tensor product of two faithful states, which is always faithful (see, e.g., the appendix of \cite{Avi82}). 
By  Lemma \ref{lem:moments}, we need only show that $u_f$ and $u_f \otimes v$ have the same $*$-moments. 
But if $w$ is any nontrivial reduced word, then $\tau_f \otimes \varphi (w(u_f \otimes v)) = \tau_f(w(u_f)) \otimes \varphi(w(v)) = 0$ because $\tau_f(w(u_f)) = 0$. 
On the other hand, clearly, $\tau_f \otimes \varphi (1 \otimes 1) = 1$. 
The result follows. 
\end{proof}

\begin{lemma}\label{lem:norm-sum-commutator}
Let $u,v$ be a pair of free Haar unitaries. Then
\[
\|u+v\|=2 \text{ and } \|vu-quv\|=2 \text{ for all $q\in\mathbb{T}$}.
\]
\end{lemma}

\begin{proof}
First, note that $u^*v$ is also a Haar unitary, so $u^*v + v^*u$ is arcsin distributed on $[-2,2]$ (see \cite[Example 1.14]{NS06}). 
Therefore, $2+u^*v + v^*u$ is arcsin distributed on $[0,4]$, from which we conclude
\[
\|u+v\|^2=\|2+u^*v + v^*u\|=4.
\]
The second claim follows from the first one together with the simple observation that $vu$ and $-quv$ are again free Haar unitaries (with respect to the same state). 
\end{proof}


\begin{lemma}\label{lem:norm-h_f}
Let $h_f:=\sum_{i=1}^d u_{f,i}+u_{f,i}^*$. 
Then we have
\[
\|h_f\|=2\sqrt{2d-1}.
\]
\end{lemma}

  This result is originally due to Kesten \cite[Theorem 3]{Kes59}, who gives a probabilistic proof; in his notation, $h_f=2d M(G,A,P)$ where $G$ is the free group $\mathbb F_d$, $A$ the canonical set of generators, $P$ assigns probability $\frac{1}{2d}$ to each generator, and $M$ is the (infinite) matrix of transition propabilities of the associated symmetric random walk on $G$. See also \cite[Theorem IV J]{AO76} for a C*-algebraic proof ($h_f$ is a \emph{free operator} in the sense \cite[Definition III B]{AO76} by \cite[Theorem III E]{AO76}) and \cite{Leh99} for a proof using free probability (the norm can be calculated from the general formula for the norm of free operators as in \cite{AO76}).

\begin{theorem}\label{thm:contr-f}
Let $a$ be a $d$-tuple of contractions on a Hilbert space $\cH$. 
Then $a\prec \frac{d}{\sqrt{2d-1}} u_f$. 
\end{theorem}

\begin{proof}
By Remark \ref{rem:contraction_unitaries}, it suffices to prove the theorem for a $d$-tuple $v\in B(\cH)^d$ of unitaries. 
In this case,  we can argue similar to \cite[Theorem 6.1]{GS20}. 
Let $u=u_f$ denote a $d$-tuple of free Haar unitaries.
First we find a state $\varphi$ on $C^*(u)$ such that $|\varphi(h_f)| = \|h_f\|$. 
By Lemma~\ref{lem:auto} we can apply a gauge automorphism and assume that $\varphi(u_i) \geq 0$ for all $i$, and hence $\varphi(h_f) = \|h_f\|$. 
Using Lemma \ref{lem:auto} again, we can replace $\varphi$ by $\frac{1}{d!}\sum_{\pi \in S_d}\varphi\circ \sigma_\pi$, which allows us to assume without loss of generality that $\varphi(u_i) = \varphi(u_1) = \alpha > 0$ for all $i$. 
Consequently, $\alpha = \frac{\|h_f\|}{2d}$ for all $i$. 
Further, by passing to the GNS representation (and recalling that $C_r^*(\bF_d)$ is simple) we may assume that $u \in B(\cH_\varphi)^d$ and that $\varphi$ is a vector state: $\varphi(a) = \langle a\xi, \xi \rangle$ for a unit vector $\xi \in \cH_\varphi$. 

Now we can form the dilation as follows. 
Put $c = \alpha^{-1} = \frac{2d}{\|h_f\|}$, and let $U_i = c u_i \otimes v_i$. 
By Lemma \ref{lem:tensor_f}, $U$ is a scalar multiple of the $d$-tuple of free Haar unitaries (one can choose an arbitrary faithful state on $C^*(v)$ to apply the lemma). 
Letting $w \colon \cH \to \cH_\varphi \otimes \cH$ be the isometry $w(h) = \xi \otimes h$,  we find that 
\[
v = c \langle u\xi, \xi\rangle v = w^* U w \prec U = c u \otimes v.
\]
By Lemma~\ref{lem:norm-h_f}, $\|h_f\|=2\sqrt{2d-1}$, and that concludes the proof.  
\end{proof} 

\begin{remark}\label{rem:contr-f}
If $\|\sum a_i+a_i^*\|=2d$, the constant in the theorem is best possible. Indeed, $a\prec cu_f$ implies $\|\sum a_i+a_i^*\|\leq c\|h_f\|$ and, thus, $c\geq\frac{2d}{\|h_f\|}=\frac{d}{\sqrt{2d-1}}$. 
In particular, the optimal dilation constant from universal unitaries to free unitaries is given by $c_{\mathrm{u},f}=\frac{d}{\sqrt{2d-1}}$, and the same goes for the optimal dilation constant from commuting unitaries to free. 
Of course, the constant is not optimal in general, as the case $a = u_f$ ($d>1$) shows. 
\end{remark}

\begin{corollary}\label{cor:c_commUni_free}
The optimal dilation constant $c_{\ru,f} = c(\ru,u_f)$ for dilating the universal tuple of unitaries to free Haar unitaries is given by
\[
c_{\mathrm{u},f}=\frac{d}{\sqrt{2d-1}}. 
\]
The optimal dilation constant $c_{0,f} = c(u_0, u_f)$ for dilating the universal commuting tuple of unitaries to free Haar unitaries is also given by the same value
\[
c_{0,f}=\frac{d}{\sqrt{2d-1}}. 
\]
\end{corollary} 

We now return to the problem that has been our primary interest: dilating unitaries to tuples of commuting unitaries. 
We begin with a lower bound for $c(u_f,u_0)$. 

\begin{theorem}\label{thm:cf0-geq-sqrt2}
For $d\geq 2$, let $c^{(d)}_{f,0} = c(u_f^{(d)}, u^{(d)}_0)$ be the dilation constant from free Haar unitaries to commuting unitaries. 
Then 
\[
c_{f,0}^{(d)}\geq 2\sqrt{1-\frac{1}{d}} .
\]
In particular, $c^{(2)}_{f,0}\geq \sqrt{2}$.
\end{theorem}

\begin{proof}
We begin with the case $d = 2$, it being considerably simpler than the general case. 
Given two unitaries $u,v$, we consider the matrix
\[
U(u,v):=\begin{pmatrix}
u^*& v\\
-v^* & u
\end{pmatrix}
\]
and find its norm to be 
\[
\left\|\begin{pmatrix}
u^* & v\\
-v^* & u
\end{pmatrix}\right\|^2
= \left\|\begin{pmatrix}
2 & [u,v]\\
[v^*,u^*] & 2
\end{pmatrix}\right\|.
\]
Therefore, $\|U(u_{0,1},u_{0,2})\|=\sqrt 2$. 
For the free Haar unitaries unitaries, note that the spectrum of
\[
\begin{pmatrix}
0 & [u_{f,1}, u_{f,2}] \\
[u_{f,2}^*, u_{f,1}^*] & 0 
\end{pmatrix}
\]
is symmetric, so using Lemma~\ref{lem:norm-sum-commutator} we get
\[
\|U(u_{f,1}, u_{f,2})\|^2 = 2+\left\|
\begin{pmatrix}
0 & [u_{f,1}, u_{f,2}] \\
[u_{f,2}^*, u_{f,1}^*] & 0 
\end{pmatrix} \right\|
= 2 + \left\|[u_{f,1}, u_{f,2}]\right\| = 4.
\]
  
To get the lower bound for $c^{(2)}_{f,0}$, suppose $u^{(2)}_{f} \prec c u^{(2)}_{0}$. 
Then $U(u_{f,1}, u_{f,2}) \prec cU(u_{0,1}, u_{0,2})$ and $2=\|U(u_{f,1}, u_{f,2})\|\leq c\|U(u_{0,1}, u_{0,2})\|=c\sqrt{2}$, so $c\geq\sqrt{2}$.

Now we consider the case $d \geq 2$. 
Given a sequence $v = (v_k)$ of unitaries on a Hilbert space $\cH$, we define operators $T_d$ on $\cH^{2^{d}} = \cH \otimes \bC^{2^d}$ recursively by setting 
\[
T_1(v):=\begin{pmatrix} 0 & v_1 \\ v_1^* & 0 \end{pmatrix}
\] 
and
\[
T_d(v) :=
\begin{pmatrix}
T_{d-1}(v) & v_d\otimes I_{2^{d-1}}\\ v_d^* \otimes I_{2^{d-1}} & -T_{d-1}(v)^*
\end{pmatrix}
\] 
(note that $T_d(v)$ is selfadjoint). 
On the one hand, plugging $v = u_0 = (u_{0,1}, \ldots, u_{0,d})$ into this construction, one inductively shows that $T_d(u_0)^*T_d(u_0)=d I_\cH \otimes I_{2^{d}}$, so $\|T_d(u_0)\|=\sqrt{d}$. 

In order to evaluate the norm of $T_d(u_f)$ obtained by plugging $v = u_f = (u_{f,1}, \ldots, u_{f,d})$ into the construction, we need to examine the structure of $T_d(v)$ a little more carefully. 
Given a sequence of unitaries $v$, we construct for all $1 \leq k \leq m$ selfadjoint unitaries $x_k^m$ recursively as follows. 
We begin by defining
\[
x_1^1 :=\begin{pmatrix} 0 & v_1 \\ v_1^* & 0 \end{pmatrix} .  
\]
Fixing $m>1$, we define for all $1\leq k < m$
\[
x_k^m :=\begin{pmatrix} x_k^{m-1} & 0 \\ 0 & -x_k^{m-1} \end{pmatrix} 
\]
and 
\[
x_m^m :=\begin{pmatrix} 0 & v_m\otimes I_{2^{d-1}} \\ v_m^*\otimes I_{2^{d-1}} & 0 \end{pmatrix} .  
\]
We now note that $T_d(v) = x^d_1 + \ldots + x^d_d$. 
Now, if we take $u_f = (u_{f,1}, \ldots, u_{f,d})$ for our sequence $v$, then using induction it can be shown that $x^d_1, \ldots, x^d_d$ are $d$ selfadjoint unitaries that have the same moments (with respect to the state $\tau_f \otimes \mathrm{tr}$) as the free Bernoulli operators, i.e., the canonical generators $\lambda(g_1), \ldots, \lambda(g_d)$ of the reduced free product group C*-algebra $C_r^*(\Gamma_d)$ where $\Gamma_d = \bZ_2 * \cdots * \bZ_2$. 
Thus, by Lemma \ref{lem:moments}, $T_d(u_f)$ is $*$-isomorphic to the operator $\sum_{i=1}^d\lambda(g_i) \in C_r^*(\Gamma_d)$, whence $\|T_d(u_f)\|=2\sqrt{d-1}$ by Lehner's formula \cite[Equation (1.1)]{Leh99} (this also follows from \cite[Theorem IV J]{AO76}, because the $\lambda(g_i)$ have the Leinert property).  
Therefore, 
\[
c_{f,0}^{(d)} \geq \frac{\|T_d(u_f)\|}{\|T_d(u_0)\|} = \frac{2\sqrt{d-1}}{\sqrt{d}} = 2\sqrt{1-\frac{1}{d}}, 
\]
and the proof is complete. 
\end{proof}

We now proceed to find an upper bound for $c_{f,0}^{(d)}$. 
\begin{theorem}\label{thm:low_freecomm}
For all $d\geq 2$, 
\[
c_{f,0}^{(d)} \leq \sqrt{2d}\frac{\sqrt{2d-1}}{d} = 2\sqrt{1-\frac{1}{2d}} .
\]
\end{theorem}

\begin{proof}
We shall use the operation of {\em polar dual}, which when applied to the matrix range of a tuple $A \in B(\cH)^d$ yields the {\em free spectrahedron} associated to $A$: 
\begin{align*}
\cW(A)^\circ 
&:= \{ X \in \bM^d : \re \sum_{j=1}^d X_j \otimes Y_j \leq 1 \,\, \textrm{ for all } Y \in \cW(A)\} \\ 
&= \{X \in \bM^d : \re \sum_{j=1}^d X_j \otimes A_j \leq 1\} =: \cD_A. 
\end{align*}
The first and third equalities above are the definitions of $\cW(A)^\circ$ and $\cD_A$, while the second equality is an easy fact \cite[Proposition 3.1]{DDSS17}. 
Now, the matrix ranges $\cW(u_f)$ and $\cW(u_0)$ are invariant under gauge actions, and therefore so are their polar duals $\cD_{u_f}$ and $\cD_{u_0}$. 
Moreover, as mentioned in Example \ref{expl:gauge}, both matrix ranges contain a neighborhood of the origin, so by \cite[Lemma 3.2]{DDSS17}, $\cW(u_f) \subseteq \cW(cu_0)$ if and only if $\cD_{cu_0} \subseteq \cD_{u_f}$. 
By \cite[Theorem 5.1]{DDSS17} (for the first equivalence below), and making use of the gauge invariance (for the third equivalence) we therefore have
\begin{align*}
u_f\prec cu_0 
&\iff 
\cW(u_f)\subseteq c \cW(u_0) \\
&\iff 
\cD_{cu_0} \subseteq \cD_{u_f} \\
&\iff 
\left\|\re \sum b_i\otimes u_{f,i}\right\|\leq c\left\|\re \sum b_i\otimes u_{0,i}\right\| \textrm{ for all } b\in \bM^d .
\end{align*}
Using Lehner's formula \cite[Corollary 1.2]{Leh99} and the operator concavity of the square-root function (see, e.g., Theorems V.1.9 and V.2.5 in \cite{BhatiaBook}), we get the following inequality 
\begin{align*}
\frac{1}{2d}\left\|\sum a_i\otimes u_{f,i} + a_i^*\otimes u_{f,i}^*\right\|
&\leq \inf_{s>0} \left\| -(1-\frac{1}{d})s+\frac{1}{2d}\sum\sqrt{s^2 + a_ia_i^*} + \sqrt{s^2 + a_i^*a_i}\right\|\\
&\leq \inf_{s>0} \left\|-(1-\frac{1}{d})s+\sqrt{s^2 +  \frac{\sum a_ia_i^* + a_i^*a_i}{2d}}\right\|\\
&=\inf_{s>0} \left\|-(1-\frac{1}{d})s+\sqrt{s^2 +  \frac{\|\sum a_ia_i^* + a_i^*a_i\|}{2d}}\right\|\\
&\leq \frac{\sqrt{2d-1}}{d}\sqrt{\frac{\|\sum a_ia_i^* + a_i^*a_i\|}{2d}}
\end{align*}
where the last inequality is obtained by taking $s=\frac{d-1}{\sqrt{2d-1}}\sqrt{\frac{\|\sum a_ia_i^* + a_i^*a_i\|}{2d}}$. (Note that in \cite[Corollary 2]{Leh97}, Lehner uses this method to prove a similar, but not quite the same inequality, so we preferred to repeat the argument.)
Therefore, 
\[
\left\|\sum a_i\otimes u_{f,i} + a_i^*\otimes u_{f,i}^*\right\|\leq \sqrt{2d}\frac{\sqrt{2d-1}}{d}\left\|\sum a_ia_i^* + a_i^*a_i\right\|^{\frac{1}{2}}.
\] 
To get the claimed inequality, note that 
\[
\sum a_ia_i^* + a_i^*a_i = \int_{z\in\mathbb T^d}\left(\sum z_ia_i + \overline{z_i}a_i^*\right)^2 \,\mathrm{d}z
\] 
and, thus,
\begin{align*}
\left\|\sum a_ia_i^* + a_i^*a_i\right\|
&\leq  \int_{z\in\mathbb T^d}\left\|\sum z_ia_i + \overline{z_i}a_i^*\right\|^2 \,\mathrm{d}z \\
&\leq \sup_{z\in\mathbb T^d}\left\|\sum z_ia_i + \overline{z_i}a_i^*\right\|^2\\
&=\left\|\sum a_i\otimes u_{0,i}+ a_i^*\otimes u_{0,i}^*\right\|^2 . 
\end{align*}
Combining everything we find that $u_f \prec cu_0$ for $c = \sqrt{2d}\frac{\sqrt{2d-1}}{d}$, and the proof is complete. 
\end{proof}

From Corollary \ref{cor:c_commUni_free} and Theorem \ref{thm:low_freecomm}, using the multiplicative triangle inequality $c(\ru, u_0) \leq c(\ru,u_f) c(u_f,u_0)$, we recover the upper bound that Passer obtained by a different method. 
\begin{corollary}[Theorem 4.4, \cite{Passer}]\label{cor:Passer}
$C_d = c(\ru^{(d)},u^{(d)}_0) \leq \sqrt{2d}$. 
\end{corollary}

We believe that the bound in the above corollary is not optimal. 
Numerical tests which were run by Matan Gibson and Ofer Israelov on random Haar unitaries suggest that $c_{f,0}^{(2)}=\sqrt{2}$, which would lead to the improved upper bound $C_2 \leq 2\sqrt{\frac{2}{3}}$ (see \cite{ShaBlog} for details).

\section{Noncommutative tori: absolute constants}\label{sec:NCT_abs} 

In this section and in the next one, we study the noncommutative tori in the context of dilation theory. 
We find bounds on the dilation constants, and then apply the general results from Section \ref{sec:metric} in combination with these bounds to recover the fact that the noncommutative tori form a continuous field of C*-algebras in a very strong sense. 

For a real and antisymmetric $d\times d$ matrix $\Theta=(\theta_{k,\ell})$, define the {\em higher dimensional noncommutative torus} $A_\Theta$ (also known as a {\em higher dimensional rotation algebra}) as the universal C*-algebra generated by $d$ unitaries $u_1, \ldots, u_d$ that satisfy the commutation relation: 
\be\label{eq:Theta_commuting}
u_\ell u_k=e^{i\theta_{k,\ell}}u_k u_\ell \,\, , \,\, k,\ell = 1, \ldots, d; 
\ee
that is, 
\[
A_\Theta:=C^*_u(u_1,\ldots, u_d:\text{unitary}, u_\ell u_k=e^{i\theta_{k,\ell}}u_k u_\ell).
\]
It is plain to see that $A_\Theta$ carries a natural gauge action as discussed in Example \ref{expl:gauge}. 

If we want to emphasize the parameters, we write $u_{\Theta,k}$ for $u_k$ and $u_\Theta = (u_{\Theta,1},\ldots, u_{\Theta,d})$. 
When $\Theta$ is equal to the $d \times d$ zero matrix $0 = 0_{d \times d}$ then $u_\Theta = u_0$ is simply the universal $d$-tuple of commuting unitaries, which we discussed in previous sections, and the notation is consistent with what we used. 

If $U = (U_1, \ldots, U_d)$ is a tuple of unitaries that satisfies \eqref{eq:Theta_commuting} (with $U_i$ instead of $u_i$, $i=k,\ell$),  then we say that $U$ {\em commutes according to $\Theta$}. 
Alternatively, we put $Q = (\exp(i\theta_{k,\ell}))$, and we say that {\em $U$ is $Q$-commuting}. 
We shall require the following lemma, to be able to deduce that certain representations of $A_\Theta$ are faithful. 

\begin{lemma}\label{lem:Qcommuting}
Let $U=(U_1,\ldots,U_d)$ be a $d$-tuple of $Q$-commuting unitaries. 
Then the following statements are equivalent. 
\begin{enumerate}[label=\textnormal{(\textit{\roman*})}]
\item\label{it:univ-univ} The canonical map $A_\Theta\to C^*(U)$ is a $*$-isomorphism.
\item\label{it:univ-aut} For all $\lambda=(\lambda_1,\ldots, \lambda_d) \in \mathbb{T}^d$ there is a $*$-automorphism $\gamma_\lambda$ of $C^*(U)$ with $\gamma_\lambda(U_k)=\lambda_kU_k$.
\item\label{it:univ-state} There is a state $\tau$ on $C^*(U)$ with $\tau(P(U))=P(0)$ for all $*$-polynomials $P$.
\end{enumerate}
\end{lemma}

\begin{proof}
Equivalence of \ref{it:univ-univ} and \ref{it:univ-state} is \cite[Lemma 4.1]{Gao18}. 
$\ref{it:univ-univ}\implies\ref{it:univ-aut}$ is obvious. 
We are left with showing $\ref{it:univ-aut}\implies\ref{it:univ-state}$. 
The map $\lambda \mapsto \gamma_\lambda(a)$ is easily seen to be continuous whenever $a$ is a $*$-polynomial in $U$, and hence by an approximation argument it is continuous for all $a \in C^*(U)$, so we can define a state
\[
\tau(a):=\int_{\mathbb{T}^d}\gamma_\lambda(a)\,\mathrm{d}\lambda.
\]
Then, clearly, $\tau(U_1^{k_1}\cdots U_d^{k_d})=0$ for $(k_1,\ldots, k_d)\neq(0\ldots 0)$.
\end{proof}

\subsection{Observations on the optimal dilation scale}\label{sec:obs-optimal}

\begin{lemma}\label{lem:trans_sym}
The optimal dilation scales
\[
c_{\Theta,\Theta'} := c(u_\Theta, u_{\Theta'})
\]
are symmetric and translation invariant. 
In particular, 
\[
c_{\Theta,\Theta'} = c_{\Theta-\Theta',0} = c_{0,\Theta-\Theta'}.
\]
\end{lemma}

\begin{proof}
We can identify $u_{\Theta + \Gamma}$ with $u_\Theta\otimes u_\Gamma$. 
Therefore every dilation of $u_\Theta$ to $cu_{\Theta'}$ gives rise to dilation of $u_{\Theta + \Gamma}$ to $cu_{\Theta' + \Gamma}$, so we get $c_{\Theta+\Gamma,\Theta'+\Gamma}\leq c_{\Theta,\Theta'}$ for all $\Theta,\Theta',\Gamma$; of course, with different choice of $\Theta,\Theta',\Gamma$, the opposite inequality also follows from this.

$A_\Theta$ and $A_{-\Theta}$ are $*$-isomorphic via $u_{\Theta,k}\mapsto u_{-\Theta,d+1-k}$. It follows that $c_{\Theta,0}=c_{-\Theta,0}$. 
By translation invariance, $c_{-\Theta,0}=c_{0,\Theta}$.
Finally,
\[
c_{\Theta,\Theta'} = c_{0,\Theta'-\Theta} = c_{\Theta-\Theta',0} = c_{\Theta'-\Theta,0} = c_{\Theta',\Theta}.
\]
\end{proof}
The above lemma allows us to concentrate on finding and formulating the values of the dilation constants 
\[
c_\Theta := c_{\Theta,0}
\]
from noncommutative to commutative, which were to begin with the constants of greatest interest. 
Recall that these constants could potentially give new information on $C_d$, as we have the lower bound $C_d \geq \sup_\Theta c_\Theta$.

\begin{lemma}
If $\alpha \geq 0$ satisfies $\alpha \leq \|\re X\|$ for all $X\in\conv(u_\Theta)$, then the linear functional $u_k\mapsto \alpha$, $1\mapsto 1$ on the unital operator space $\operatorname{span}(1,u_\Theta)$ is contractive. 
\end{lemma}

\begin{proof}
For brevity, let us write $u = u_\Theta$. 
Let $b + \sum a_ku_k\in\operatorname{span}(1,u)$ and assume without loss of generality that $\sum |a_k|=1$. 
Under the condition $\alpha \leq \|\re X\|$ for all $X\in\conv(u_\Theta)$, we get (using the existence of the gauge automorphisms $\gamma_\lambda$) 
  \begin{align*}
    \left\|b + \sum a_ku_k\right\|
    &=\left\||b| + \sum |a_k|u_k\right\|\\
    &\geq\left\|\re\left(|b| + \sum |a_k|u_k\right)\right\|\\
    &=|b| + \left\|\re\sum |a_k| u_k\right\|\\
    &\geq|b|+ \alpha \geq \left|b +\sum a_k \alpha \right|.
  \end{align*}
\end{proof}

\begin{theorem}\label{thm:cTheta}
For all $\Theta$, we have
\[
c_\Theta = \frac{1}{\inf\{\|\re X\|:X\in\conv(u_\Theta)\}}.
\]
\end{theorem}

\begin{proof}
By Lemma \ref{lem:trans_sym}, $c_{\Theta}=c_{0,\Theta}$, so we obtain the lower bound by considering the latter constant. 
Suppose that $u_0 \prec c u_\Theta$. 
For every $X = \sum t_k u_{\Theta,k} \in \conv(u_\Theta)$ we have that $Y := \sum t_k u_{0,k} \prec cX$ and therefore $\re Y \prec c\re X$. 
It is easy to see that $\|\re Y\|= \sup_{z\in \bT^d} | \re \sum_k t_k z_k | = 1$. 
On the other hand,  $\|\re Y\| \leq c \|\re X\|$, or $c\geq \frac{\|\re Y\|}{\|\re  X\|}=\frac{1}{\|\re X\|}$, so $c_\Theta$ is bounded below as claimed: 
\[
c_\Theta = c_{0,\Theta} \geq \frac{1}{\inf\{\|\re X\|:X\in\conv(u_\Theta)\}}.
\] 

To show that the above inequality is an equality, we will construct a commuting normal dilation consisting of unitaries of norm $\frac{1}{\alpha}$, where $\alpha := \inf\{\|\re X\|:X\in\conv(u_\Theta)\}$. 
The construction is similar to the proof of Theorem \ref{thm:contr-f}. 
We note that $u_\Theta \otimes u_{-\Theta} := (u_{\Theta,1} \otimes u_{-\Theta,1},\ldots,  u_{\Theta,d} \otimes u_{-\Theta,d})$ is a tuple of commuting unitaries. 
By the previous lemma, and recalling that $A_\Theta$ and $A_{-\Theta}$ are $*$-isomorphic, we find that there is a state $\varphi$ on $A_{-\Theta}$ with $\varphi(u_{-\Theta,k}) = \alpha$ for all $k$. 
If we define $U = \frac{1}{\alpha} u_\Theta\otimes u_{-\Theta}$, then $U$ is $\frac{1}{\alpha}$ times a tuple of commuting unitaries, and ${\id} \otimes \varphi (U) = u_\Theta$. 
From this is follows that $u_\Theta$ has the required dilation, so that $c_\Theta \leq \frac{1}{\alpha}$, as required. 
\end{proof}

One may reformulate the above by defining 
\[
x_{\Theta,k} =\frac{1}{2}\left( u_{\Theta,k} + u_{\Theta,k}^*\right) \quad, \quad k =1, \ldots, d, 
\]
and then noting that elements of the form $\re X$ for $X \in \conv(u_\Theta)$ are precisely the elements in $\conv(x_\Theta)$. 
Thus, we seek the convex combination $\sum t_k x_{\Theta,k}$ of minimal norm. 

In \cite{GS20}, the case $d = 2$ was studied. 
With $u_\theta$ denoting the universal pair satisfying $u_{\theta,2}u_{\theta,1} = e^{i\theta} u_{\theta,1} u_{\theta,2}$, the main result (Theorem 6.3) in that paper can be stated as follows: 
\[
c_\theta = \frac{4}{\|u_{\theta,1}+u_{\theta,1}^*+u_{\theta,2}+u_{\theta,2}^*\|}, 
\]
in other words, the element of minimal norm in $\conv(x_\theta)$ is the uniform mixture $\frac{1}{4}h_\theta$, where $h_\theta:= 2 x_{\theta,1} + 2 x_{\theta,2} = u_{\theta,1}+u_{\theta,1}^*+u_{\theta,2}+u_{\theta,2}^*$. 
This result was used in \cite[Section 6]{GS20} to give the best currently known lower bound for $C_2$: 
\[
C_2 \geq \max_\theta c_\theta \geq c_{\theta_s} \approx 1.5437772, 
\]
where $\theta_s = 2\pi(\sqrt{2}-1)$ is the value of $\theta$ where we conjecture that the maximum is attained. 

The following result gives a crude upper bound for $c_\Theta$ (the fact that this bound is far from optimal can be seen using the ideas going into Corollary \ref{cor:lower_bound}). 

\begin{proposition}\label{prop:tensor}
Let $\Theta = (\theta_{k,\ell})_{k,\ell=1}^d$ be a real antisymmetric $d \times d$ matrix. 
Then, 
\[
c_\Theta \leq \prod_{\ell = 2}^d \max_{1 \leq k \leq \ell-1} c_{\theta_{k,\ell}}  .
\]
\end{proposition}
\begin{proof}
By the proof of \cite[Theorem 4.2]{Gao18}, a continuous family of universal representations of $A_\Theta$ for all $\Theta$ can be defined recursively as follows. 
Let $\hat{\Theta}$ denote the $(d-1)\times (d-1)$ matrix $\hat{\Theta} = (\theta_{k,\ell})_{k,\ell=1}^{d-1}$ obtained by removing the last row and column of $\Theta$. 
If the representation $U^{(d-1)}(\hat{\Theta}) = (U_1^{(d-1)}(\hat{\Theta}), \ldots, U_{d-1}^{(d-1)}(\hat{\Theta}))$ of $u_{\hat{\Theta}}$ has already been defined on $\cK$, then the representation $U^{(d)}(\Theta) = (U_1^{(d)}(\Theta), \ldots, U_d^{(d)}(\Theta))$ of $u_\Theta$ is defined on $\cK \otimes \cH \otimes \cdots \otimes \cH$ by 
\begin{align*}
U_1^{(d)}(\Theta) &= U_1^{(d-1)}(\hat{\Theta}) \otimes u(\theta_{1,d}) \otimes I \otimes \cdots \otimes I \\
U_2^{(d)}(\Theta) &= U_2^{(d-1)}(\hat{\Theta}) \otimes I \otimes u(\theta_{2,d}) \otimes \cdots \otimes I \\
 &\,\,\,\vdots  \\
U_{d-1}^{(d)}(\Theta) &= U_{d-1}^{(d-1)}(\hat{\Theta}) \otimes I \otimes I \otimes \cdots \otimes u(\theta_{d-1,d}) ,
\end{align*}
and
\[
\,\,\,\,\,\, U_d^{(d)}(\Theta) = I \otimes v(\theta_{1,d}) \otimes v(\theta_{2,d}) \otimes \cdots \otimes v(\theta_{d-1,d}) , 
\]
where $\theta \mapsto u(\theta)$ and $\theta \mapsto v(\theta)$ are the H{\"o}lder continuous paths of universal $e^{i\theta}$-commuting unitaries on $\cH$ that were constructed in \cite{HR95}. 
This family of representations is continuous in $\Theta$, by induction. 
The fact that $U(\Theta)$ commutes according to $\Theta$ is plain, and the fact that it is universal follows from simplicity of $A_\Theta$ for a dense set of $\Theta$s (see \cite{Gao18} for further explanations). 

Now, if $V = (V_1, \ldots, V_{d-1})$ is a commuting family of unitaries and $c \geq 1$ is such that $U^{(d-1)} \prec c V$, and $(w(\theta), z(\theta))$ is a pair of commuting unitaries such that $(u(\theta), v(\theta)) \prec c_\theta (w(\theta), z(\theta))$, then we define a tuple $N$ of commuting normals on $\cK \otimes \cH \otimes \cdots \otimes \cH$ by 
\begin{align*}
N_1 &= V_1 \otimes w(\theta_{1,d}) \otimes I \otimes \cdots \otimes I \\
N_2 &= V_2 \otimes I \otimes w(\theta_{2,d}) \otimes \cdots \otimes I \\
 &\,\,\,\vdots  \\
N_{d-1} &= V_{d-1} \otimes I \otimes I \otimes \cdots \otimes w(\theta_{d-1,d}) ,
\end{align*}
and
\[
N_d = I \otimes z(\theta_{1,d}) \otimes \cdots \otimes z(\theta_{d-1,d}) . 
\]
We see that $N$ is a tuple of commuting normal contractions, and that $U^{(d)}(\Theta) \prec C N$ for $C = c \max_{1\leq k \leq d-1}c_{\theta_{k,d}}$. 
\end{proof}

\subsection{The dilation constant for $q$-commuting $d$-tuples}\label{sec:q_commuting_d}

Our next proposition gives a precise value of $C_\Theta$ for matrices $\Theta$ with a a constant value above the diagonal, i.e., for tuples $u = (u_1, \ldots, u_d)$ for which there is a $q \in \bT$ such that $u_\ell u_k = q u_k u_\ell$ for all $k<\ell$. 
\begin{proposition}\label{prop:sametheta}
Let $\Theta$ be a real antisymmetric $d \times d$ matrix such that $\theta_{k,\ell} = \theta$ for all $k<\ell$. 
Then
\[
c_\Theta = \frac{2d}{\|u_{\Theta,1}+u_{\Theta,1}^* + \ldots + u_{\Theta,d} + u_{\Theta,d}^*\|}.
\]
\end{proposition}
\begin{proof}
Let us write $u = u_\Theta$, and define $h = u_1 + u_1^* + \ldots + u_d + u_d^*$. 
By Theorem \ref{thm:cTheta}, we know that $c_\Theta = \frac{1}{\alpha}$, where $\alpha = \inf\{\|\re X\| : X \in \conv(u)\}$. 
Therefore, $c_\Theta \geq 2d/\|h\|$. 
Moreover, by the second half of the proof of that theorem, we know how to construct a commuting normal dilation of norm $c$ once we find a state $\varphi$ on $A_\Theta$ such that $\varphi(u_k) = 1/c$ for all $k$. 
Therefore, our task boils down to finding a state $\varphi$ such that $\varphi(u_k) = \frac{\|h\|}{2d}$ for all $k$; this will show that $c_\Theta \leq \frac{2d}{\|h\|}$, and the proof will be complete. 

Let $\psi$ be a state such that $|\psi(h)| = \|h\|$. 
By gauge invariance we may assume that $\psi(u_k) \geq 0$, and so also $\psi(u_k) = \psi(u_k^*)$, for all $k$. 
For the special kind of $\Theta$ we are considering, the map $u_k \mapsto u_{k+1}$ for $k=1, \ldots, d-1$
and $u_d \mapsto u_1^*$ extends to an automorphism $\sigma$ of $A_\Theta$ that fixes $h$. 
Letting $\varphi = \frac{1}{d}\sum_{k=0}^{d-1} \psi \circ \sigma^k$ we find the desired state that satisfies $\varphi(u_k) = \frac{\|h\|}{2d}$ for all $k$. 
\end{proof}

From the previous proposition we obtain a new lower bound for $C_3 = c(\ru,u_0)$, improving the previously known bound $C_3 \geq \sqrt{3} \approx 1.732$ given by \eqref{eq:Cd}. 
\begin{corollary}\label{cor:lower_bound}
Let $\Sigma$ be the unique $3\times 3$ antisymmetric matrix that has $1$s above the diagonal, so that for every $\theta \in \bR$, $\theta \Sigma$ is the antisymmetric matrix with $\theta$ above the diagonal. 
Then 
\[
C_3 \geq \max_\theta c_{\theta \Sigma} \geq c_{\frac{6\pi}{7}\Sigma} \geq 1.858 . 
\]
\end{corollary}
\begin{proof}
The first and second inequalities are obvious. 
We guessed that the angle $\theta = 2\pi\times \frac{3}{7}$ will give a relatively large value of $c_{\theta \Sigma}$ by running some numerical simulations. 
To obtain a reliable lower bound $c_{\frac{6\pi}{7}\Sigma} \geq 1.858$ we combined numerics and theory as follows. 

Let us write $q = e^{i\theta}$, where $\theta = 2\pi \frac{m}{n}$, where $\gcd(m,n) = 1$. 
Define
\[
X := \begin{pmatrix} 
q & & &  & \\ 
 & q^2 & & & \\ 
&  & q^3 & & \\ 
 & & & \ddots & \\ 
 & & & & 1 \end{pmatrix},
\qquad
 Y := 
 \begin{pmatrix}
 & 1 & & & \\
 & & 1 & & \\ 
 & & & \ddots  & \\
 & & & & 1 \\
 1 & & & & 
 \end{pmatrix} .
 \]
For brevity, let us write $u = u_{\theta\Sigma}$. 
One can show that every irreducible representation of $A_{\theta \Sigma}$ is determined by $u_1 \mapsto \alpha X$, $u_2 \mapsto \beta XY$ and $u_3 \mapsto \gamma Y$, where $\alpha, \beta, \gamma \in \bT$ (see Appendix~\ref{sec:repr-A(3)q}). 
Therefore, to compute $\|u_1+u_1^*+u_2+u_2^*+u_3+u_3^*\|$ we have to solve the finite dimensional optimization problem
\[
\max\{\|\alpha X + (\alpha X)^* + \beta XY + (\beta XY)^* + \gamma Y + (\gamma Y)^*\| : \alpha, \beta, \gamma \in \bT\}. 
\]
It is not a difficult matter to numerically find this maximum to reasonable precision, and one can then estimate an upper bound on the resulting error in the estimate for $c_\theta$. 
The lower bound $c_{\frac{6\pi}{7}} \geq 1.858$, together with the simple MATLAB code and the analysis leading to it can be found in \cite{GPSS_script}. 
\end{proof}

\section{Noncommutative tori: continuity of the dilation constants}\label{sec:NCT_cont}

\subsection{Representing noncommutative tori by Weyl unitaries}\label{sec:repr-high-dimens}

We shall make use of the symmetric Fock space and the Weyl unitaries (see \cite{Par12}; be aware that we use inner products which are linear in the first argument). 
For a Hilbert space $H$ let
\[
\Gamma(H):=\bigoplus_{k=0}^{\infty} H^{\otimes_s k}
\]
be the symmetric Fock space over $H$.
The exponential vectors $e(x):=\sum_{k=0}^{\infty}\frac{1}{\sqrt{k!}} x^{\otimes k}, x\in H$ form a linearly independent and total subset of $\Gamma(H)$.
Clearly, $\langle e(x),e(y)\rangle = e^{\langle x,y\rangle}$ for all $x,y\in H$. For $z\in H$ we define the Weyl unitary $W(z)\in B(\Gamma(H))$ which is determined by
\[W(z) e(x)=e(z+x) \exp\left(-\frac{\|z\|^2}{2} - \langle x,z\rangle \right)\]
for all exponential vectors $e(x)$.
A simple calculation shows that $W(z),W(y)$ commute up to the phase factor $e^{2i \im \langle y,z\rangle}$; to be precise:
\[
W(y) W(z) = e^{2i \im \langle y,z\rangle} W(z) W(y).
\]

\begin{proposition}\label{prop:Weyl_univ}
Let $x_1,\ldots, x_d\in H$ be linearly independent and $2\im\langle x_\ell,x_k\rangle = \theta_{k,\ell}$. Then the C*-algebra generated by $W(x_1),\ldots,W(x_d)$ is $*$-isomorphic to $A_\Theta$. 
\end{proposition} 

\begin{proof}
Note that the unitaries $W(x_1), \ldots, W(x_d)$ commute according to $\Theta$. 
By Lemma \ref{lem:Qcommuting}, it is enough to show the existence of gauge automorphisms $\gamma_\lambda$ as in condition \ref{it:univ-aut} of the lemma. 
Let $\lambda_k=e^{it_k}$. Because $x_1, \ldots, x_d$ are linearly independent, there exists a vector $x\in H$ with $\langle x_k,x\rangle = \frac{i}{2}t_k$ for all $k$. Therefore, $W(x)^*W(x_k)W(x)=e^{it_k}W(x_k)$ and therefore conjugation with $W(x)$ is the gauge automorphism $\gamma_\lambda$ on the C*-algebra generated by $W(x_1), \ldots, W(x_d)$.
\end{proof}

\begin{lemma}
Let $\Theta=(\theta_{k,\ell})$ be a real and antisymmetric $d\times d$ matrix. Then there exists a Hilbert space $H$, $\dim H=d$ and linearly independent $x_1,\ldots, x_d\in H$ such that
\[
2\im\langle x_\ell,x_k\rangle = \theta_{k,\ell}.
\]
\end{lemma}

\begin{proof}
Let $H$ be a Hilbert space with an orthonormal basis $e_1,\ldots, e_d$. 
We define $x_k$ recursively:
\begin{itemize}
\item $x_1:=e_1$
\item $x_{k}:=e_{k}+ \widetilde{x}_{k}$ with $\widetilde{x}_{k}\in\operatorname{span}\{e_1,\ldots, e_{k-1}\}$, $\langle x_\ell, \widetilde{x}_{k}\rangle=\frac{i}{2}\theta_{k,\ell}$
\end{itemize}
to get linearly independent $x_1,\ldots,x_d$ with   $2\im\langle x_\ell,x_k\rangle = \theta_{k,\ell}$.
\end{proof}

\begin{lemma}
Let $\Theta=(\theta_{k,\ell})$, $\Theta'=(\theta'_{k,\ell})$ be real and antisymmetric $d\times d$ matrices. Then there exist Hilbert spaces $H\subset K$, $\dim H=d$, $\dim K = 2d$ and linearly independent vectors $z_1,\ldots, z_d\in K$, such that $x_k:=P_H z_k$ ($k=1, \ldots, d$) are linearly independent, and such that with $y_k:=P_{H^{\perp}}z_k$ the following conditions are satisfied:
\begin{enumerate}[label=\textnormal{(\arabic*)}]
\item $2\im\langle z_\ell,z_k\rangle = \theta_{k,\ell}'$
\item $2\im\langle x_\ell,x_k\rangle = \theta_{k,\ell}$
\item $\|y_k\|^2=\frac{1}{2}\|\Theta'-\Theta\|$
\end{enumerate}
\end{lemma}

\begin{proof}
Construct $x_1,\ldots,x_d$ as in the previous lemma.
Next, note that $\frac{i}{2}(\Theta'-\Theta)$ is selfadjoint, so $\frac{1}{2}\|\Theta'-\Theta\|I_d + \frac{i}{2}(\Theta'-\Theta)\geq0$. 
Thus, there exists a $d \times d$ matrix $Y$ such that
\[
Y^*Y=\frac{1}{2}\|\Theta'-\Theta\|I_d + \frac{i}{2}(\Theta'-\Theta).
\] 
Let $y_k$ be the $k$th column of $Y$, so that $\|y_k\|^2=\frac{1}{2}\|\Theta'-\Theta\|$ and $2\im\langle y_\ell,y_k\rangle = \theta'_{k,\ell}-\theta_{k,\ell}$.
Finally, define $K = H \oplus \bC^d$ and put $z_k:=x_k\oplus y_k \in K$. 
Then the $z_k$ are clearly linearly independent, as they project onto the linearly independent $x_k$, and
\[
2\im\langle z_\ell,z_k\rangle = 2\im\langle x_\ell,x_k\rangle + 2\im\langle y_\ell,y_k\rangle = \theta_{k,\ell} + (\theta'_{k,\ell}-\theta_{k,\ell})=\theta'_{k,\ell}, 
\]
as required. 
\end{proof}

We now reach the main result of this section. 
\begin{theorem}\label{thm:Weyl-rep}
In the notation of the previous lemma:
\begin{enumerate}[label=\textnormal{(\arabic*)}]
\item $A_\Theta$ has a representation by Weyl unitaries, $u_k\mapsto W(x_k)$;
\item $A_{\Theta'}$ has a representation by Weyl unitaries, $u_k\mapsto W(z_k)$;
\item $W(x_1),\ldots, W(x_d)$ is a compression of $cW(z_1),\ldots cW(z_d)$ with $c=e^{\frac{1}{4}\|\Theta-\Theta'\|}$. 
\end{enumerate}
Consequently, 
\[
c_{\Theta,\Theta'} \leq e^{\frac{1}{4}\|\Theta-\Theta'\|}. 
\]
\end{theorem}

\begin{proof}
Let $p$ denote the projection of $K$ onto $H$, and let $p^\perp$ denote the projection onto the orthogonal complement $H^\perp$. 
Consider the symmetric Fock spaces $\Gamma(H) \subset \Gamma(K)$ with $P$ the
projection onto $\Gamma(H)$. 
By Proposition \ref{prop:Weyl_univ}, the C*-algebras generated by $W(x_1), \ldots, W(x_d)$ and $W(z_1), \ldots, W(z_d)$ are $*$-isomorphic to $A_\Theta$ and $A_{\Theta'}$, respectively. 

Note that for exponential vectors we have
$Pe(z)=e(pz)$.
A calculation then shows that for every $k,\ell = 1, \ldots, d$, the following conditions are satisfied: 
\begin{enumerate}[label=\textnormal{(\arabic*)}]
\item $W(z_k), W(z_\ell)$ commute up to the phase factor $e^{2i \im \langle z_\ell, z_k \rangle}$;
\item $PW(z_k)\big|_{\Gamma(H)}= e^{-\frac{\|p^\perp  z_k\|^2}{2}}W(pz) = e^{-\frac{\|\Theta-\Theta'\|}{4}}W(x_k)$;
\item $W(x_k)$ and $W(x_\ell)$ commute up to the phase factor $e^{2i \im \langle x_\ell, x_k\rangle}$.
\end{enumerate}
This completes the proof. 
\end{proof}

\subsection{Continuous representations of noncommutative tori}\label{sec:cont_reps}

It is known that the Weyl operators do not form a norm continuous family, in the sense that $\|W(x) - W(y)\| \geq \sqrt{2}$ if $x \neq y$ (see \cite[Proposition 2.2]{PetzBook}). 
On the other hand, Haagerup and R{\o}rdam showed that there is a $\frac{1}{2}$-H{\"o}lder continuous path $\theta \mapsto u_\theta = (u_{\theta,1}, u_{\theta,2})$ such that $u_{\theta,2}u_{\theta,1} = e^{i\theta}u_{\theta,1} u_{\theta,2}$ for all $\theta$ \cite{HR95} (the analogous result for the higher dimensional rotation algebras follows from this, see \cite[Theorem 1.1]{Gao18}). 
A key step in Haagerup and R{\o}rdam's result was to show that there is a constant $K>0$ such that 
\begin{align}
\label{eq:HRest}
\distance_{\mathrm{HR}}(u_\theta,u_{\theta'}) \leq K |\theta - \theta'|^{1/2},
\end{align}
see \cite[Theorem 4.9]{HR95}. 
They proved this inequality by showing first that unbounded selfadjoint operators satisfying the Heisenberg commutation relation can be boundedly approximated by commuting unbounded selfadjoint operators. 
From \eqref{eq:HRest} (essentially), Haagerup and R{\o}rdam then proceed to construct the H{\"o}lder continuous path $\theta \mapsto u_\theta$. 
Our methods provide an alternative and direct proof of the higher dimensional analogue of \eqref{eq:HRest}. 

\begin{corollary}\label{cor:nctori_cont}
For every two real antisymmetric matrices $\Theta, \Theta'$, 
\[
\distance_{\mathrm{HR}}(u_\Theta,u_{\Theta'}) \leq 28 \|\Theta - \Theta'\|^{1/2}. 
\]
\end{corollary}
\begin{proof}
By Theorem \ref{thm:Weyl-rep}, $\distance_{\mathrm{D}}(u_\Theta, u_{\Theta'}) = \log c_{\Theta,\Theta'} \leq \frac{1}{4}\|\Theta - \Theta'\|$, so we invoke Theorem \ref{thm:dilandHR} and find $\distance_{\mathrm{HR}}(u_\Theta,u_{\Theta'})\leq 56\distance_{\mathrm{D}}(u_\Theta, u_{\Theta'})^{1/2}\leq 28  \|\Theta - \Theta'\|^{1/2}$.
\end{proof}

The above corollary is powerful enough to prove quite easily that for every selfadjoint $*$-polynomial, the spectrum $\sigma(p(u_\Theta))$ is a $\frac{1}{2}$-H{\"o}lder continuous map, and from this one easily obtains the known results that the higher dimensional rotation algebras $A_\Theta$ form a continuous field. 
We will not elaborate here on this method (see \cite[Section 4]{GS20} for details in the case $d=2$), since much stronger conclusions can be drawn by invoking an ingenious technique from \cite{HR95}. 

In \cite{Gao18} (following \cite{HR95}) it was proved that there exists a norm continuous map $\Theta \mapsto U(\Theta) \in B(\cH)^d$, such that $U(\Theta) \sim u_\Theta$ and
\[
\|U(\Theta) - U({\Theta'})\| \leq K \|\Theta - \Theta'\|^{1/2} .
\]
Corollary \ref{cor:nctori_cont} brings us very close to recovering that result. 
To derive the existence of a norm continuous embedding $\Theta \mapsto U(\Theta)$ from the above corollary, we first observe, following \cite{Gao18}, that the construction in the proof of Proposition \ref{prop:tensor} allows us to reduce to the case of $d=2$. 
To prove the case $d=2$, we note that Corollary \ref{cor:nctori_cont} implies that the map $\theta \mapsto u_\theta = (u_{\theta,1}, u_{\theta,2})$ projects onto a $\frac{1}{2}$-H{\"o}lder continuous path in $\cU(2)$. 
In \cite[Section 5]{HR95}, Haagerup and R{\o}rdam showed that this particular path can be lifted to a $\frac{1}{2}$-H{\"o}lder continuous path into $\cU(\cH) \times \cU(\cH)$ for some separable Hilbert space $\cH$. 
In the following appendix, we will modify Haagerup and R{\o}rdam's techniques to prove that for all $\alpha \in (0,1)$, every $\alpha$-H{\"o}lder continuous path into $\cU(d)$ can be lifted to an $\alpha$-H{\"o}lder continuous path into $\cU(\cH)^d$ for some separable Hilbert space $\cH$.

\appendix
\section{Lifting H{\"o}lder continuous paths}
\label{sec:lift-lip_-paths}

In this section, we track down the proofs of \cite[Lemmas 5.2 -- 5.4]{HR95} to prove that H{\"o}lder continuous paths into $\cU(d)$ can be lifted to H{\"o}lder continuous paths (of the same exponent) into $\cU(\cH)^d$. 
Most of the ideas that follow are from \cite{HR95}, but the results there are not in a form directly applicable, so we record the necessary modified proofs here. 
The following remarkable lemma is key.

\begin{lemma}[{\cite[Lemma 5.1]{HR95}}]\label{HR-lemma}
Let $M\subset B(H)$ be a von Neumann algebra with properly infinite
commutant $M’$, and let $u\in M$ be unitary. Then there is a smooth path $u(t),
t \in [0, 1]$, of unitaries in $B(H)$, such that
\begin{enumerate}[label=\textnormal{(\roman*)}]
\item $u(0)=1$ and $u(1)= u$;
\item $\|u'(t)\| \leq 9$;
\item $\|[u(t), a]\| \leq 4\|[u, a]\|$;
\item $\|[u'(t), a]\| \leq 9\|[u, a]\|$;
\item $\left\|\frac{\mathrm d}{\mathrm d t} u(t)au(t)^*\right\| \leq 45 \|[u, a]\|$;
\end{enumerate}
for all $t\in [0, 1]$ and all $a\in M$.  
\end{lemma}

Recall that $\cU(d)$ is the quotient of $\cU(\cH)^d$ by the equivalence relation that identifies $*$-isomorphic tuples (see Section \ref{sec:metric})
and that two $d$-tuples of unitaries $U$ and $V$ are equivalent 
if and only if the infinite ampliations of $U$ and $V$ are approximately unitarily equivalent
. 
Note that a unitary tuple $U$ is the infinite ampliation of some unitary tuple if and only if it is unitarily equivalent to its own second ampliation. In this case, for every unitary $W$, $WUW^*$ is again an infinite ampliation. This allows us to assume that all unitary tuples we construct will be infinite ampliations, in particular, two tuples will be equivalent if and only if they are approximately unitarily equivalent.

\begin{lemma}\label{lem:lift-sequence}
  Let $u_0,\ldots, u_k\in \cU(d)$ and $\varepsilon >0$. Then there are $U_0,\ldots, U_k$ with $U_i\sim u_i$ and $\|U_{i+1}-U_i\|\leq \distance_{\mathrm{HR}}(u_{i+1},u_i) + \varepsilon$ for all $i=0,\ldots, k-1$. 
\end{lemma}

\begin{proof}
  For $k=0$ there is nothing to prove. Now assume we already have such a realization for $u_0,\ldots, u_{k-1}$. As a simple consequence of the definition of $\distance_{\mathrm{HR}}$, we can choose $V_{k-1}\sim u_{k-1}$ and $V_k\sim u_k$ with $\|V_{k}-V_{k-1}\|\leq \distance_{\mathrm{HR}}(u_{k},u_{k-1}) + \frac{\varepsilon}{2}$. Take the infinite ampliations of all $U_i$ and $V_i$ if necessary to assure that $U_{k-1}$ and $V_{k-1}$ are approximately unitarily equivalent. Then there is a unitary $W$ such that $\|U_{k-1}-WV_{k-1}W^*\|<\frac{\varepsilon}{2}$. Therefore, putting $U_k = WV_k W^*$, we have 
\begin{align*}
\|U_k-U_{k-1}\|&\leq \|WV_kW^*-WV_{k-1}W^*\| +  \|WV_{k-1}W^*-U_{k-1}\| \\ 
&\leq \distance_{\mathrm{HR}}(u_k,u_{k-1})+\varepsilon
\end{align*}
as required. 
\end{proof}

\begin{lemma}\label{lem:lift-sequence-start-end-fixed}
Let $u_0,\ldots, u_k\in \cU(d)$ with $\distance_{\mathrm{HR}}(u_i,u_{i+1})<\delta$ for some $\delta>0$. 
Fix two representatives $U_0\sim u_0$, $U_k\sim u_k$ acting on the same Hilbert space $\cH$, such that the commutant of $U_0,U_k$ is properly infinite. 
Then there are unitaries $U_1,\ldots, U_{k-1}$ on $\cH$ such that
\[
\|U_{i+1}-U_i\|\leq 227\delta + 45\frac{\|U_k-U_0\|}{k}
\]
and $U_i\sim u_i$ for all $i=0,\ldots, k-1$, and such that the commutant of $U_0,U_1,\ldots, U_k$ is properly infinite. 
\end{lemma}

\begin{proof}
We closely follow the proof of \cite[Lemma 5.2]{HR95}. 
Since the commutant of $U_0, U_k$ is properly infinite, we can write $\cH=H_1\otimes H_2\otimes H_3$ with infinite dimensional $H_1,H_2,H_3$ so that $U_0,U_k\in B(H_1)^d\otimes 1 \otimes 1 = (B(H_1)\otimes \mathbb C1 \otimes \mathbb C1)^d$. 
 We first apply Lemma~\ref{lem:lift-sequence} to obtain $\widetilde U_0,\ldots, \widetilde U_k$ with $\widetilde U_i\sim u_i$ and $\|\widetilde U_{i+1}-\widetilde U_i\| < \delta$. 
Without loss of generality we assume that the $\widetilde U_i$ are infinitely ampliated and belong to $B(H_1)^d\otimes 1 \otimes 1$. 
The unitaries $U_0$ and $\widetilde U_0$ have the form $X \otimes 1 \otimes 1$ and $Y \otimes 1 \otimes 1$, where $X\otimes 1$ and $Y\otimes 1$ are approximately unitarily equivalent, so  there is a unitary $V\in B(H_1\otimes H_2) \otimes 1$ such that $\|U_0-V\widetilde U_0V^*\| < \delta$. 
Put $\overline U_i:=V\widetilde U_i V^*$. 
By the triangle inequality,
\[
\|\overline U_i-\overline U_j\| \leq k\delta
\]
for all $i,j\in\{0,\ldots, k\}$. 
With $\Delta:=\|U_k-U_0\|$, it follows that
\[
\|\overline U_k-U_k\| \leq (k+1)\delta + \Delta. 
\]
By the same reasoning as for $U_0$, there is a unitary $W\in B(H_1\otimes H_2) \otimes 1$ such that $\|U_k-W\overline U_kW^*\|<\delta$. 
Accordingly,
\begin{align*}
\|\overline U_k-W\overline U_kW^*\| 
&\leq \|\overline U_k-\overline U_0\|+\|\overline U_0-U_0\|+\|U_0-U_k\|+\|U_k-W\overline U_kW^*\| \\
& \leq (k+2)\delta + \Delta
\end{align*}
and, for every $i=0,\ldots, k$, 

\begin{align*}
\|\overline U_i-W\overline U_iW^*\|
&\leq \|\overline U_i-\overline U_k\|+ \|\overline U_k-W\overline U_kW^*\|+\|W\overline U_kW^*-W\overline U_iW^*\| \\
& \leq (3k+2)\delta + \Delta. 
\end{align*}
Now apply Lemma~\ref{HR-lemma} with $M = W^*(\overline U_0, \ldots, \overline U_k)$ and $u = W$ to obtain a smooth path $W(t)$, $t\in[0,1]$ of unitaries in $B(H_1\otimes H_2)\otimes 1$ such that $W(0)=1, W(1)=W$ and
\[
\left\|\frac{\mathrm d}{\mathrm d t} W(t)\overline U_iW(t)^*\right\|\leq 45\|[W,\overline U_i]\| \leq 45((3k+2)\delta + \Delta).
\]
Put $\hat U_i:=W(i/k)\overline U_iW(i/k)^*$ for $i=0,\ldots,k$. Then
\begin{align*}
\|\hat U_{i+1}-\hat U_{i}\|&\leq \|W({\textstyle\frac{i+1}{k}})(\overline U_{i+1}-\overline U_i)W({\textstyle\frac{i+1}{k}})^*\|+\left\|\int_{\frac{i}{k}}^{\frac{i+1}{k}}\frac{\mathrm{d}}{\mathrm{d}t}(W(t)\overline U_iW(t)^*) \mathrm d t\right\|\\
    &\leq \delta + \frac{45((3k+2)\delta + \Delta)}{k}\leq 226\delta + 45\frac{\Delta}{k}.
\end{align*}
Put $U_i:=\hat U_i$ for $i=1,\ldots,k-1$. 
From $\|U_0-\hat U_0\| < \delta$ and $\|U_k-\hat U_k\| < \delta$ it follows that
\[
\|U_{i+1}-U_i\|\leq 227\delta + 45\frac{\Delta}{k}.
\]
By construction, all the unitary tuples $U_0, \ldots, U_k$ are infinite ampliations and contained in $B(H_1\otimes H_2)^d \otimes \mathbb C1$, and so their commutant is properly infinite.  
\end{proof}

Let $\alpha \in (0,1]$. 
Recall that a path $\gamma\colon [a, b]\to X$ from an interval $[a,b]$ to a metric space $(X,d)$ is said to be {\em H{\"o}lder continuous with exponent $\alpha$}, or simply {\em $\alpha$-H{\"o}lder}, if there is a constant $C$ such that
\be\label{eq:Lipa}
d(\gamma(s),\gamma(t))\leq C|t-s|^\alpha
\ee
for all $s,t\in[a,b]$. 
H{\"o}lder continuous paths with exponent $\alpha = 1$ are simply called {\em Lipschitz continuous}.

\begin{lemma}\label{Lem:Gamma}
Fix $k\in\mathbb N$, $k>1$ and consider the set $\Gamma=\bigcup_{n=0}^{\infty}\Gamma_n\subset[0,1]$ with
\[
\Gamma_n=\left\{\frac{j}{k^n}~\middle\vert~ j=0,1\ldots,k^n\right\}.
\] 
Let $X$ be a complete metric space and $\gamma\colon \Gamma\to X$ a function such that there are constants $C_1>0$ (independent of $n$) and $0<\alpha\leq 1$ with
\[
d(\gamma(s),\gamma(t))\leq C_1 |t-s|^{\alpha}
\]
for all $s,t\in \Gamma_n$ with $t-s=k^{-n}$ (i.e., for every pair of adjacent points in $\Gamma_n$). 
Then $\gamma$ has an $\alpha$-H{\"o}lder continuous extension $\overline{\gamma}$ to $[0,1]$ such that 
\[
d(\overline\gamma(s),\overline\gamma(t))\leq C|t-s|^{\alpha}
\]
for some $C>0$ and all $s,t\in[0,1]$.
\end{lemma}

\begin{proof}
The case of $\alpha = 1$ is elementary, and in fact a Lipschitz continuous extension exists with $C = C_1$. 
We treat the case $0 < \alpha < 1$. 

First, note that for all $s,t\in\Gamma_n$ with $|t-s|=\ell k^{-n}$ and $0\leq\ell<k$, the triangle inequality yields 
\be\label{eq:dgamma}
d(\gamma(t),\gamma(s))\leq \ell C_1 k^{-n\alpha}\leq k^{1-\alpha}C_1 (\ell k^{-n})^{\alpha}=C_2|t-s|^{\alpha}
\ee
with $C_2:=k^{1-\alpha}C_1$.

For arbitrary $s<t\in\Gamma$, there is a smallest $M$ with $s,t\in\Gamma_{M}$ and a unique $N$ such that $k^{-N}\leq |t-s|< k^{-(N-1)}$ (let us not worry about the case $s = 0, t = 1$). 
Recursively, one finds intermediate points $s_i,t_i\in \Gamma_{i}$, $i\in\{N,N+1,\ldots, M\}$ such that
\[
s=s_{M}\leq s_{M-1}\leq \ldots \leq s_{N}\leq t_{N} \leq\cdots \leq t_{M}=t,
\]
$s_{i-1}-s_{i}=a_i k^{-i}$, $t_{i}-t_{i-1}=b_i k^{-i}$ with $a_i,b_i\in\{0,\ldots, k-1\}$; clearly, also $t_{N}-s_{N}=ck^{-N}$ with $c\in\{0,\ldots, k-1\}$, because we assumed that $|t-s|< k^{-(N-1)}$.
The triangle inequality together with \eqref{eq:dgamma} now yields
\begin{align*}
d(\gamma(t),\gamma(s))
& \leq C_2 \left(\sum_{i=N+1}^M (a_ik^{-i})^{\alpha}+ \sum_{i=N+1}^M (b_ik^{-i})^{\alpha} +(ck^{-N})^{\alpha}\right)\\
& \leq C_2 \left(\sum_{i=N+1}^M (k^{-i+1})^{\alpha}+ \sum_{i=N+1}^M (k^{-i+1})^{\alpha} +(k^{-N+1})^{\alpha}\right)\\
& \leq 2C_2 k^\alpha \sum_{i=N}^{M} k^{-i\alpha} \\
& \leq 2C_1 k\left(\sum_{i=0}^{\infty} (k^{-\alpha})^{i}\right) k^{-N\alpha} \\
& \leq \frac{2kC_1}{1-k^{-\alpha}} |t-s|^{\alpha}
\end{align*}
and so $d(\gamma(t),\gamma(s)) \leq C|t-s|^{\alpha}$ with $C=\frac{2kC_1}{1-k^{-\alpha}}$. 
We have shown that $\gamma\colon \Gamma\to X$ is uniformly $\alpha$-H{\"o}lder continuous, and from this the statement of the lemma follows easily. 
\end{proof}

\begin{theorem}\label{thm:Holder_ext}
Let $0 < \alpha < 1$, let $\gamma\colon[0,1]\to (\cU(d),\distance_{\mathrm{HR}})$ be an $\alpha$-H{\"o}lder path, and suppose that $U_0\sim \gamma(0)$ and $U_1 \sim \gamma(1)$ are two representatives with a properly infinite commutant. 
Then there is an $\alpha$-H{\"o}lder path $U\colon [0,1]\to \cU(\cH)^d$ with $U(0) = U_0$,\, $U(1) = U_1$ and $U(t) \sim \gamma(t)$ for all $t \in [0,1]$. 
\end{theorem}

\begin{proof}
Let us assume that $C$ is such that a strict inequality holds in \eqref{eq:Lipa} for all $s \neq t$. 
We may and do assume that $C \geq 1$. 
Fix a $k\in\mathbb N$ with $k\geq 90^{\frac{1}{1-\alpha}}$ and consider the set $\Gamma$ as defined in Lemma~\ref{Lem:Gamma}. 
Define $U(0) = U_0$ and $U(1) = U_1$. 
We can successively apply Lemma~\ref{lem:lift-sequence-start-end-fixed} to construct $U(\frac{j}{k^{n+1}}) \sim \gamma(\frac{j}{k^{n+1}})$ (with properly infinite commutant) for the points $\frac{j}{k^{n+1}}$ that belong to $\Gamma_{n+1}$ but not to $\Gamma_n$. 
Let $d_n:=\max_j \|U(\frac{j+1}{k^n})-U(\frac{j}{k^n})\|$. 
Then
\begin{align}
\left\|\textstyle U\left(\frac{j+1}{k^{n+1}}\right)-U\left(\frac{j}{k^{n+1}}\right)\right\|\leq A\frac{C}{k^{(n+1)\alpha}} + B\frac{d_n}{k},\label{eq:1}
\end{align}
for $A = 227$ and $B = 45$, because
\[
\distance_{\mathrm{HR}}\left({\textstyle \gamma\left(\frac{j}{k^{n+1}}\right),\gamma\left(\frac{j+1}{k^{n+1}}\right)}\right)
< 
\frac{C}{k^{(n+1)\alpha}}=:\delta
\]
for all $j$.
It follows inductively that
\be\label{eq:2}
d_n\leq \frac{2AC}{k^{n\alpha}} .
\ee
Indeed, $d_0 \leq 2 \leq 2B C$, and
\begin{align*}
d_{n+1} 
&\leq A\frac{C}{k^{(n+1)\alpha}}+B \frac{d_n}{k} \\
&\leq \frac{AC}{k^{(n+1)\alpha}} + \frac{B}{k}\left(  \frac{2AC}{k^{n\alpha}} \right) \\ 
&\leq \left(A+ B\frac{2A}{k^{1-\alpha}}\right)\frac{C}{k^{(n+1)\alpha}}
\end{align*}
which gives the result because we chose $k\geq 90^{\frac{1}{1-\alpha}}$, so that $\frac{B}{k^{1-\alpha}} = \frac{45}{k^{1-\alpha}} \leq \frac{1}{2}$. 
This way a function $\Gamma\to \cU(\cH)^d$ is defined.
Finally, we find that for adjacent $s,t\in\Gamma_n$  
\[\|U(t)-U(s)\|\leq d_{n}\leq \frac{2AC}{k^{n\alpha}} = 2AC|t-s|^{\alpha},\]
so an application of Lemma~\ref{Lem:Gamma} finishes the proof. 

Now, $U(t) \sim \gamma(t)$ for all $t \in \Gamma$ by construction. 
By continuity of $\gamma$ and $U$, we find that $\distance_{\mathrm{HR}}(U(t), \gamma(t)) = 0$  for all $t \in [0,1]$. 
Using that $\distance_{\mathrm{HR}}$ is a metric (Proposition \ref{prop:metrics}), we conclude that $U(t) \sim \gamma(t)$. 
\end{proof}

It is natural to ask whether every Lipschitz continuous path $\gamma\colon[0,1]\to (\cU(d),\distance_{\mathrm{HR}})$ can be lifted to a Lipschitz continuous path into $\cU(\cH)^d$. 
The following theorem gives a partial answer to this question. 

\begin{theorem}\label{thm:Lip_ext}
Let $\gamma\colon[0,1]\to (\cU(d),\distance_{\mathrm{HR}})$ be a Lipschitz continuous path, and suppose that $U_0\sim \gamma(0)$ and $U_1 \sim \gamma(1)$ are two representatives with a properly infinite commutant. 
Then there is a continuous path $U\colon [0,1]\to \cU(\cH)^d$ that is $\alpha$-H{\"o}lder for all $0<\alpha<1$, satisfying $U(0) = U_0$,\, $U(1) = U_1$ and $U(t) \sim \gamma(t)$ for all $t \in [0,1]$. 
\end{theorem}

\begin{proof}
The proof is modelled on that of Theorem \ref{thm:Holder_ext}. 
Again, we may choose $C\geq 1$ such that strict inequality holds in \eqref{eq:Lipa} for all $s \neq t$ and $\alpha = 1$. 
Let $k\in \bN$, whose value will be chosen later. 
As above, we define $U(0) = U_0$ and $U(1) = U_1$, however now we will refine the sequence of partitions more rapidly. 
Let $k_n = k^{n^2}$. 
We apply Lemma~\ref{lem:lift-sequence-start-end-fixed} successively to construct $U(\frac{j}{k_{n+1}})$ for the points $\frac{j}{k_{n+1}}$ that belong to $\Gamma_{(n+1)^2}$ but not to $\Gamma_{n^2}$. 
Let us set $d_n:=\max_j \|U(\frac{j+1}{k_{n}})-U(\frac{j}{k_{n}})\|$; the estimation of these quantities will now be slightly more delicate. 
We have that 
\begin{align}
\left\|\textstyle U\left(\frac{j+1}{k_{n+1}}\right)-U\left(\frac{j}{k_{n+1}}\right)\right\|
\leq ACk^{-(n+1)^2} + B\frac{d_n}{k^{2n+1}}, \label{eq:3}
\end{align}
because $\gamma$ is Lipschitz, and because $\Gamma_{(n+1)^2}$ is obtained from $\Gamma_{n^2}$ by partitioning every interval into $\frac{k_{n+1}}{k_n} = \frac{k^{(n+1)^2}}{k^{n^2}} = k^{2n+1}$ sub-intervals. 
We will now prove that $d_n \leq D \left(\frac{1}{k_n}\right)^{-\alpha_n}$, where $D \geq 2$ and $\alpha_n \nearrow 1$ is a sequence that will be soon determined. 
Start from $d_0 \leq 2 \leq D$. 
Proceeding inductively, we plug $d_n \leq D \left(\frac{1}{k_n}\right)^{-\alpha_n}$ in \eqref{eq:3} we find 
\begin{align*}
d_{n+1} 
&\leq A C k^{-(n+1)^2} + B Dk^{-n^2\alpha_n - 2n -1} \\
&= \left(AC + BDk^{n^2(1-\alpha_n)}\right)k^{-(n+1)^2} \\
&= \left(AC + BDk^{n^2(1-\alpha_n)}\right) k^{-(n+1)^2(1-\alpha_{n+1})}\left(\frac{1}{k_{n+1}}\right)^{\alpha_{n+1}} .
\end{align*}
Our inductive step will be complete if we can arrange that 
\[
\left(AC + BDk^{n^2(1-\alpha_n)}\right) k^{-(n+1)^2(1-\alpha_{n+1})} \leq D . 
\]
It suffices to find the parameters that satisfy
\[
k^{n^2(1-\alpha_n)-(n+1)^2(1-\alpha_{n+1})} \leq \frac{D}{AC + BD}. 
\]
First, we choose $D$ large enough so that $\frac{D}{AC + BD} > \frac{1}{2B}$. 
Next, we note that if we choose $\alpha_n = 1 - n^{-1/2}$, then 
\[
n^2(1-\alpha_n)-(n+1)^2(1-\alpha_{n+1}) = n^{3/2} - (n+1)^{3/2}\leq -1 
\]
for all $n$. 
Finally, we choose $k$ large enough so that $k^{-1} < \frac{1}{2B}$. 
With these parameters in place, we see that the inductive step holds. 

We have therefore defined a map 
\[
U \colon \Gamma \to \cU(\cH)^d
\]
such that for every pair of adjacent points $s,t \in \Gamma_{n^2}$, 
\[
\|U(s) - U(t)\| \leq D |s-t|^{\alpha_n} , 
\]
with $\alpha_n = 1-n^{-1/2}$. 
Note that $D$ does not depend on $n$. 
To show that $U$ extends to a path $U \colon [0,1] \to \cU(\cH)^d$ that is H{\"o}lder continuous for all exponents $<1$, we will need to modify the proof of Lemma \ref{Lem:Gamma}.

What we are going to prove is that $U$ extends to a path $U \colon [0,1] \to \cU(\cH)^d$ that is locally $\alpha$-H{\"o}lder continuous for every $0 < \alpha <1$, with a uniform (local) constant. 
This will be done by showing that $U$ itself satisfies a H{\"o}lder continuity condition on intervals in $\Gamma$. 
It suffices to check for $|s-t| \leq k^{-5}$, so that in the following considerations we consider $\Gamma_{n^2}$ only for $n \geq 5$, which guarantees $2n + 1 - n^{3/2} < 0$.

The difference from Lemma \ref{Lem:Gamma} is that now in the passage from $\Gamma_{n^2}$ to $\Gamma_{(n+1)^2}$ we add $k^{2n+1}-1$ points to every interval. 
Now, for $n \geq 5$, and for all distinct $s,t\in\Gamma_{n^2}$ with $|t-s|=\ell k^{-n^2}$ and $0 < \ell < k^{2n+1}$, we write $\beta_n = \alpha_n - n^{-1/2} = 1 - 2n^{-1/2}$, and by the triangle inequality
\begin{align}\label{eq:betanineq}
\|U(t) - U(s) \| \leq \ell D k^{-n^2\alpha_n} = \ell D k^{-n^{3/2}} k^{-n^2 \beta_n} \leq  D (\ell k^{-n^2})^{\beta_n} \leq  D |t-s|^{\beta_n}
\end{align}
because $\ell k^{-n^{3/2}} \leq k^{2n + 1 - n^{3/2}} \leq 1 \leq \ell^{\beta_n}$.

For arbitrary $s<t\in\Gamma$, there is a smallest $M$ with $s,t\in\Gamma_{M^2}$ and a unique $N$ such that $k^{-N^2}\leq |t-s|< k^{-(N-1)^2}$. 
Recursively, one finds intermediate points $s_m, t_m \in \Gamma_{m^2}$, $m \in \{N,N+1,\ldots, M\}$ such that
\[
s=s_{M}\leq s_{M-1}\leq \ldots \leq s_{N}\leq t_{N} \leq\cdots \leq t_{M} = t,
\]
$s_{m-1}-s_{m} = a_m k^{-m^2}$, $t_{m}-t_{m-1} = b_m k^{-m^2}$ with $a_m, b_m \in \{0,\ldots, k^{2m+1}-1\}$. 
Since we assumed that $|t-s|< k^{-(N-1)^2}$, we also have $s_{N}-t_{N} = ck^{-N^2}$ with $c\in\{0,\ldots, k^{2N+1} - 1\}$.
Writing $\delta_n = \beta_n - n^{-1/2}$, using \eqref{eq:betanineq} we obtain
\begin{align*}
d(\gamma(t),\gamma(s))
& \leq D \left(\sum_{m=N+1}^M (a_m k^{-m^2})^{\beta_m}+ \sum_{m=N+1}^M (b_m k^{-m^2})^{\beta_m} +(ck^{-N^2})^{\beta_N}\right)\\
& \leq D \left(\sum_{m=N+1}^M (k^{-m^2 +2m +1})^{\beta_m}+ \sum_{m=N+1}^M (k^{-m^2 + 2m + 1})^{\beta_m} +(k^{-N^2 + 2N + 1})^{\beta_N}\right)\\
& \leq D \left(\sum_{m=N+1}^M (k^{-m^2})^{\delta_m}+ \sum_{m=N+1}^M (k^{-m^2})^{\delta_m} +(k^{-N^2})^{\delta_N}\right)\\
& \leq 2 D \sum_{m=N}^{M} k^{-m^2 \delta_N} \leq 2 D \left(\sum_{m=0}^{\infty} k^{-m^2 \delta_N} \right) k^{-N^2\delta_N} \leq K |t-s|^{\delta_N} , 
\end{align*}
with $K:= 2 D \left(\sum_{m=0}^{\infty} k^{-m^2 \delta_5} \right)$, which does not depend on $N$. Thus, we get 
\[
\| U(t) - U(s) \| \leq K|t-s|^{\delta_N}
\]
for all $s,t \in \Gamma$ with $|s-t| < k^{-(N-1)^2}$ with $\delta_N \nearrow 1$. 

Letting $N = 5$, we see that on the intersection of $\Gamma$ with an interval of length $k^{-16}$ the map $U$ is $\delta_{5}$-H{\"o}lder. 
It follows that $U$ extends to a uniformly continuous path $U \colon [0,1] \to \cU(\cH)^d$, and using our estimate for general $N\geq 5$, it follows that $U$ is $\delta_{N}$-H{\"o}lder with constant $K$ on every open interval of length $k^{-(N-1)^2}$. 

It follows that $U$ is $\alpha$-H{\"o}lder for every $0<\alpha<1$. 
The assertion $U(t) \sim \gamma(t)$ for all $t \in [0,1]$ follows as in the proof of Theorem \ref{thm:Holder_ext}. 
\end{proof}

\section{Representations of $A^{(3)}_q$}
\label{sec:repr-A(3)q}

In this section, we prove the claim made in the proof of Corollary \ref{cor:lower_bound} concerning irreducible representations of \(A^{(3)}_q:=A_{\theta\Sigma}=C_{\mathrm u}^*(u_1,u_2,u_3: \text{unitary}, u_\ell u_k=q u_k u_\ell \text{ for $k<\ell$})\):

\begin{proposition}
  Let $q=e^{2\pi i\frac{m}{n}}$ be a primitive $n$-th root of unity. Then every irreducible representation of $A^{(3)}_q$ is of the form $u_1\mapsto \alpha X$, $u_2\mapsto \beta Z=\beta XY, u_3\mapsto \gamma Y$ where $\alpha,\beta,\gamma\in \mathbb T$ and
  \[X = \begin{pmatrix} 
q & & &  & \\ 
 & q^2 & & & \\ 
&  & q^3 & & \\ 
 & & & \ddots & \\ 
 & & & & 1 \end{pmatrix} ,\ 
 Y = 
 \begin{pmatrix}
 & 1 & & & \\
 & & 1 & & \\ 
 & & & \ddots  & \\
 & & & & 1 \\
 1 & & & & 
\end{pmatrix},\ 
Z= XY = 
 \begin{pmatrix}
 & q & & & \\
 & & q^2 & & \\ 
 & & & \ddots  & \\
 & & & & q^{n-1} \\
 1 & & & & 
 \end{pmatrix}.\]
\end{proposition}

\begin{proof}
  A representation of $A^{(3)}_q$ is determined by three unitaries $U_1, U_2, U_3$ such that $u_i \mapsto U_i$ for $i=1,2,3$. 
We note that $U_2 U_1^n = q^n U_1^n U_2$, etc., and we find that $U_i^n$ commutes with $U_j$ for all $i,j$. 
Therefore, if $U$ determines an irreducible representation, then $U_i^n$ is a scalar multiple of the identity for all $i$. 
It then follows that $C^*(U)$ is of dimension less than or equal to $n^3$. 
By irreducibility, again, it must act on a finite dimensional space $\cK$ (of dimension $\dim \cK \leq n^{3/2}$).

Given such a representation, $u_1 \mapsto U_1$ and $u_3 \mapsto U_3$ determines a representation of $A^{(2)}_q$ on a finite dimensional space, and so is a direct sum of irreducibles; the irreducible summands are known to be of the form $u_1\mapsto \alpha_i X,u_3\mapsto \gamma_i Y$ (cf.\ \cite{BocBook}) and we conclude that $\cK$ breaks up as $\cK = \bC^n \oplus \cdots \oplus \bC^n=(\mathbb C^n)^m$, and $U_1, U_3$ break up as block diagonal operators $U_1 = \bigoplus_{i=1}^m \alpha_i X$ and $U_3 = \bigoplus_{i=1}^m \gamma_i Y$ (where the $\alpha$s and $\gamma$s are on the unit circle). 
The operator $U_2$ has a block decomposition $U_2 = (W_{ij})_{i,j=1}^m$ where every $W_{ij}$ is an $n \times n$ matrix. 
We have the relations $U_\ell U_k = q U_k U_\ell$ for all $k<\ell$. 
For $k=1, \ell=2$, this gives 
\[
W_{ij}\alpha_j X = q \alpha_i X W_{ij}. 
\]
This implies that $\lambda = q\alpha_i/\alpha_j$ must be an $n$th root of unity (e.g., by checking that $\lambda^{-1}$ is an eigenvalue of $X$ corresponding to the vector $W_{ij} v$, where $v$ is an eigenvcetor of $X$ corresponding to $1$). 
It follows that $\alpha_i/\alpha_j$ is an $n$th root of unity. 
Therefore, there exists some $\alpha \in \bT$ such that $\alpha_j = q^{m_j} \alpha$ for all $j$. 
Repeating these considerations with $U_3$ replacing $U_1$, we find that there is some $\gamma\in\mathbb T$ such that $\gamma_j = q^{n_j} \gamma$ for all $j$.
Since the pair $(X,Y)$ is unitarily equivalent to the pair $(q^s X, q^t Y)$ for $s,t\in\mathbb Z$, we may conjugate by a block unitary to pass to $U_1 =  \bigoplus_{i=1}^m \alpha X = \alpha X \otimes I_m$ and $U_3 =  \bigoplus_{i=1}^m\gamma Y = \gamma Y \otimes I_m$. 

After this simplification, we return to the relations $U_\ell U_k = q U_k U_\ell$ for all $k<\ell$. 
For $k=1, \ell=2$, we now find
\begin{equation}
\label{eq:XW}
W_{ij} X = q X W_{ij},
\end{equation}
and for $k=2,\ell=3$ we have 
\begin{equation}
\label{eq:YW}
YW_{ij} = q W_{ij} Y. 
\end{equation}
Let us fix $1\leq i,j \leq m$ and write $W_{ij} =:W = (w_{st})_{s,t=1}^n$. 
From \eqref{eq:XW}, we find that for all $1 \leq s,t \leq n$, 
\[
w_{st} q^t = q^{s+1}w_{st}. 
\]
This means that $w_{st} = 0$ unless $t = s+1 \mod n$. 
So $W$ has the same structure as $Y$. 
Now from \eqref{eq:YW} we find
\[
w_{s+1,t} = qw_{s,t-1}, 
\]
where addition in the indices is carried out ${\operatorname{mod} n}$. 
Therefore, $W$ is a scalar multiple of
 \[
Z = \begin{pmatrix}
 & q & & & \\
 & & q^2 & & \\ 
 & & & \ddots  & \\
 & & & & q^{n-1} \\
 1 & & & & 
 \end{pmatrix} = X Y.
 \]
 We conclude that $U_2 = Z \otimes B$ with $B$ an $m \times m$ matrix (which must be unitary because $U_2$ is unitary). 
For every subspace $\cL \subseteq \mathbb C^m$ that is invariant for $B$, the space $\cH \otimes \cL$ is invariant for the representation determined by $U$, thus $U$ is irreducible if and only if $m = 1$, so that $U_2=Z\otimes B=\beta Z$ with $\beta\in\mathbb T$.
\end{proof}

\section*{Acknowledgments}

The authors would like to thank Benjamin Passer for some useful feedback.


 \linespread{1.25}


\linespread{1}
\setlength{\parindent}{0pt}

\providecommand{\bysame}{\leavevmode\hbox to3em{\hrulefill}\thinspace}
\providecommand{\MR}{\relax\ifhmode\unskip\space\fi MR }
\providecommand{\MRhref}[2]{%
  \href{http://www.ams.org/mathscinet-getitem?mr=#1}{#2}
}
\providecommand{\href}[2]{#2}

\end{document}